
\magnification=\magstep1 
\baselineskip=12pt 

\def\ind {{\lim \limits_{\longrightarrow}}} 
\def \pro {{\lim\limits_{\longleftarrow}}} 
\input amstex 
\documentstyle{amsppt} 
\overfullrule=0pt 
\vsize=500pt 

\document

\def\g{{\gamma}} 
\def\a{{\alpha}} 
\def\b{{\beta}}

\def\gen{{\frak{g}}}

\def\lra{\,\longrightarrow\,}

\def\a{{\alpha}} 
\def\o{{\omega}} 
\def\l{{\lambda}} 
\def\b{{\beta}} 
\def\g{{\gamma}} 
 
\def\eps{{\varepsilon}}

\def\1b{{\bold 1}}

\def\Ob{{\bold O}}

\def\Obj{{\roman{Ob}}}

\def\Oplus{{\ts\bigoplus}}

\def\J{{\roman J}}

\def\gr{\roman{gr}} 
 
\def\Spec{\roman{Spec\,}}

\def\Db{{\bold{D}}}

\def\sqr{{\!\sqrt{\,\,\,}}}

\def\Schb{{\bold{Sch}}} 
\def\Fschb{{\bold{Fsch}}} 
\def\Gr{{\roman{Gr}}}

\def\cGr{{{\text {co-}}\Gr}} 
\def\Ischb{{\bold{Isch}}}

\def\Lrsb{{\bold{Lrs}}} 

\def\Fsetb{{\bold{Fset}}} 
\def\Setb{{\bold{Set}}} 
\def\Aut{\roman{Aut}}

\def\max{{\text{max}}} 
\def\Hom{\roman{Hom}} 
\def\Homu{\underline{\roman{Hom}}}
\def\Der{\roman{Der}}

\def\Id{\roman{Id}} 
 
\def\CDOc{\Cc\Dc\Oc} 
\def\Detc{\Dc et} 
\def\CDOb{{\bold{CDO}}} 
\def\Zar{{\roman{Zar}}} 
\def\Et{{\roman{\acute{Et}}}}

\def\End{\text{End}\,} 
\def\Quot{\roman{Quot}}

\def\Ker{\text{Ker}\,}

\def\max{\text{max}\,} 
\def\Gr{{\text{Gr}}}

\def\AA{{\Bbb A}} 
 
\def\CC{{\Bbb C}}

\def\GG{{\Bbb G}}

\def\ZZ{{\Bbb Z}}

\def\Cc{{\Cal C}} 
\def\Dc{{\Cal D}} 
\def\Ec{{\Cal E}} 
\def\Fc{{\Cal F}} 
\def\Gc{{\Cal G}} 
\def\Hc{{\Cal H}} 
 
\def\Kc{{\Cal K}} 
\def\Lc{{\Cal L}} 
\def\Mc{{\Cal M}} 
\def\Nc{{\Cal N}} 
\def\Oc{{\Cal O}}

\def\Vc{{\Cal V}} 
\def\Wc{{\Cal W}} 
\def\Xc{{\Cal X}} 
\def\Yc{{\Cal Y}} 
\def\Zc{{\Cal Z}}
\def\Homc{{\Hc om}}

\def\and{{\quad\text{and}\quad}}

\def\ts{\textstyle}

\def\qed{\hfill $\sqcap \hskip-6.5pt \sqcup$}

\def\Wedge{{\ts\bigwedge}}

\def\lra{{{\longrightarrow}}}

\def\u1{{\underline 1}}

\centerline{\bf FORMAL LOOPS IV : CHIRAL DIFFERENTIAL OPERATORS} 

\vskip .5cm 

\centerline {\bf M. Kapranov, E. Vasserot} 

\vskip 2cm 

\noindent {\bf (0.1)}  The goal of this paper, which builds on earlier papers [KV1-3],
is to relate the theory of sheaves of chiral differential operators (CDO) [GMS1-2]
with the determinantal anomaly on the loop space.

Sheaves of CDO studied in {\it loc. cit.} are curvilinear versions of the fundamental vertex algebra
called the Heisenberg module [K]. To be precise, $A_N$, the loop Heisenberg algebra
of $\CC^N$, is generated by elements (creation and annihilation operators)
$$ a^i_n, b^i_n, \quad i=1, ..., N, \,\, n\in \ZZ,$$
subject to the relations
$[a^i_n, b^j_m] = \delta_{i,j}\delta_{m, -n}$, while other commutators are zero. 
The Heisenberg module $V_N$ is the cyclic (vacuum) $A_N$-module $A_N\bigl/ A_n\{ b^i_{<0}, a^i_{\leq 0}\}$.  
 
It is known that $V_N$ is a vertex algebra and the question is to find analogs of $V_N$
when $\CC^N$ is replaced by an arbitrary $N$-dimensional complex algebraic manifold $X$.
Intuitively, the generators of $A_N$ are the coefficients of $\CC^N$-valued ``fields"
$$a(t) = \bigl(a^1(t), ..., a^N(t)\bigr) , \quad b(t)= 
\bigl(b^1(t), ..., b^N(t)\bigr), \quad
a^i(t) = \sum _n a^i_n t^n \,\,\, {\text {etc}.}
$$
while in a sheaf of CDO one would like to deal with $X$-valued fields.
In other words, we want to construct a sheaf of vertex algebras on $X$ whose
``local model" is $V_N$.

\vskip .3cm

\noindent {\bf (0.2)} The problem of constructing sheaves of CDO was analyzed in
[GMS1-2]  from the change of variables point of view, i.e.,  by constructing an action
of the group of coordinate changes in $\CC^N$ on (some completion of) $V_N$.
This action turns out to be a projective one which leads to an anomaly, or obstruction
which, in general, prevents one from gluing such a sheaf globally. 
 To be precise, sheaves of CDO form a gerbe $\Cc\Dc\Oc_X$ on $X$
with lien $\Omega^{2, cl}_X$, the sheaf of closed 2-forms, and the class of this gerbe,
i.e., the obstruction to the existence of a global object, is
$$[\Cc\Dc\Oc_X] = {1\over 2} c_1(X)^2 - c_2(X)\,\,\, \in \,\,\, H^2(X, \Omega^{2, cl}_X).
\leqno (0.2.1)$$

\vskip .3cm

\noindent {\bf (0.3)} In the present paper we relate $\Cc\Dc\Oc_X$ to $\Lc X$, the space of formal
loops on $X$ constructed in [KV1]. This is an ind-scheme which is an algebro-geometric
model for the space of free loops $LX = C^\infty(S^1, X)$. The tangent spaces of $\Lc X$, while
infinite-dimensional, possess a ``Tate structure", see [D], and thus we have the determinantal gerbe
$\Dc et_{\Lc X}$ on $\Lc X$ with lien $\Oc^*_{\Lc X}$. Our main result, Theorem 3.4.2, says that
$\Dc et_{\Lc X}$ and $\Cc\Dc\Oc_X$ are in a sense, identified. The identification is achieved by
means of the {\it symplectic action homomorphism} 
$$S: \Omega^{2, cl}_X \to \Oc^*_{\Lc X}, \quad \omega\mapsto \exp\left( \int d^{-1}(\omega)\right),$$
see [KV3], (2.4.3). If $\omega$ is a symplectic form, then $S(\omega)$, as a function on the
(formal) loop space is the exponential of what is usually called the symplectic action functional
corresponding to $\omega$. The image of $S$ consists of functions on $\Lc X$ which are
``factorizable", see {\it loc. cit} . Factorization is a crucial concept in the approach to vertex algebras developed
by Beilinson and Drinfeld [BD1]. In our context focusing on the subsheaf of factorizable functions
(which makes $S$ to be an isomorphism) corresponds to the requirement that our curvilinear versions of $V_N$
be sheaves of vertex algebras.

\vskip .3cm

\noindent {\bf (0.4)} The main change of variables formula of [MSV], [GMS1] which led to (0.2.1), is somehow
analogous to the famous Ito formula in the theory of stochastic integrals [McK]: in both cases
it is the second derivatives of the transformation functions which make their way into the place where
one normally expects first derivatives only. In fact, the fields $a(t), b(t)$ are remindful of the
Brownian motion, so the problem of constructing sheaves of CDO is analogous to that of constructing
the Brownian motion on a curved manifold [E]. One cannot help noticing  further similarities
between the normal ordering used in [GMS1-2] to regularize the meaning of the $X$-valued $a(t), b(t)$
and the anticipatory choice of the middle point in the Riemann sums approximating the Ito integrals [McK]. 
In both cases a different choice of regularization (say, the symmetric one)
is a priori possible but it would change the answer
in a not so essential way. 

\vskip .3cm

\noindent {\bf (0.5)} Our main technical point, which was already mentioned in [KV1], is this.
The would-be sheaf of chiral differential operators on $X$ is to be considered as the pushdown
(space of sections)
of the $\Dc$-module $\delta_{\Lc^0 X}$ on $\Lc X$ formed by delta-functions supported on
the subscheme $\Lc^0 X$ of Taylor loops.  However, while $\delta_{\Lc ^0 X}$ is well defined
as an object of a certain formal category $\Db$, the concept of the space of global sections
for objects of this category is not immediately defined at all. In our approach, it is this
concept of global sections which has to be regularized. The reason is the behavior of
$\Dc$-modules under the operations of direct and inverse image: say, for  a closed smooth embedding
$i: Z\hookrightarrow W$ and a right $\Dc$-module $\Mc$ on $Z$  there is no embedding of
sections of $\Mc$ into those of its direct image $i_!\Mc$: there is a correction factor
in the form of the determinant of the normal bundle $\Nc_{Z/W}$. In the infinite-dimensional situation
when we take inductive limits of such embeddings, the resulting ambiguities are
neatly described by  the determinantal gerbe $\Dc et_{\Lc X}$. 

Our first step is thus to  associate to any (local) object $\Ec$ of $\Dc et_{\Lc X}$ the functor of global
sections $\Gamma_\Ec$ from $\Db$ to vector spaces. This is done in Section 1 for  $\Dc$-modules on
a class of ind-schemes 
including the $\Lc X$ which we call locally locally compact ({\it sic!}). This class is related to 
the so-called reasonable ind-schemes of [BD2], see also [D]. The first ``locally" in the term means locally with
respect to the Zariski topology while the ``locally compact" part afterwards refers to an algebraic
counterpart of local compactness for ind-pro-objects, as formulated by Kato, see [KV1]. 
Thus Section  1 is devoted to the study of determinantal anomaly (of the contangent bundle)
for smooth ind-schemes of this class. 

\vskip .3cm

\noindent {\bf (0.6)} Our second main step is establishing the conditions for the space $\Gamma_\Ec(\delta_{\Lc^0 X})$
to be a vertex algebra. Here, as in [KV1], we use the formulation of the theory by means of
factorization structures, developed in [BD1] but we have to investigate the new level of
coherence for factorization: that of { gerbes}. This required rewriting the axioms of [BD1] in
a new way, which is suitable for such next level questions. This is done in Section 2.
We present a natural polyhedral axiomatics which can, in principle, be pushed as far down  the
coherence route as desired: we have a natural system of polyhedra in all dimensions,
similarly to the Stasheff polyhedra for associativity.  

Thus, factorization properties are present in our theory on all levels: the ind-scheme
$\Lc X$ gives rise to a factorization semigroup, the determinantal gerbe is then a factorization
gerbe in our new sense, and we can  speak about its {\it factorizing objects}. The $\Gamma$-functors
 corresponding to such objects give vertex algebras. Further, $\Hom$-sheaves between factorizing
 objects are  what we call factorizing line bundles, and among sections of factorizing line bundles there
 is a distinguished class of factorizing sections. In the particular case of the trivial
 factorizing line bundle the sections (i.e., invertible functions) that are factorizing, are precisely
 the functions in the image of $S$ from (0.3). 
 
 Finally, in Section 3 we put the two threads together and establish the relation between the
 stack formed by factorizing objects of the determinantal gerbe and the stack (gerbe) of CDO,
 proving our main theorem.

 \vskip .3cm
 
 \noindent {\bf  (0.7)}  The first author would like to acknowledge support from the NSF, Universit\'e Paris-7 and the
 Max-Planck-Institut f\"ur Mathematik in Bonn.

\vfill\eject

\centerline{\bf 1. Sections of D-modules over ind-schemes.}

\vskip 1cm

\subhead{\bf (1.1) Locally locally compact ind-schemes}\endsubhead

We work over the field $\CC$ of complex numbers. 
All rings will be assumed to contain 
$\CC$ and all schemes will be schemes over $\CC$. 
Let $\Schb$ be the category of schemes and $\Fschb$ 
be the full subcategory of schemes of finite type. 

We will follow the same conventions about ind-schemes as in \cite{KV1,2}. 
In addition, in this paper we will assume that the indexing filtering
poset of an ind-scheme is at most countable. 
Let $\Ischb$ be the category of  such ind-schemes.
As in \cite{KV2, sect.~2}, to any ind-scheme $Y$ we associate its Zariski site
$Y_\Zar$ and a sheaf of pro-algebras $\Oc_Y$ on $Y_\Zar$. Objects
of $Y_\Zar$ will be called open sub-ind-schemes in $Y$. Thus if
$Y= ``\ind"Y^\a$ an open sub-ind-scheme in $Y$ is the same
as a compatible collection of open sub-schemes $U^\a\subset Y^\a$.

\proclaim{(1.1.1) Definition} 
An ind-scheme $Y$ is called : 

(a) 
discrete if it can be represented as 
$Y=``\ind" Y^\a$ where 
$(Y^\a)$ is a filtering inductive system over $\Fschb$ 
such that each map $i^{\a\a'}\,:\,Y^{\a}\to Y^{\a'}$, 
$\a\leq\a'$, is a closed embedding, 

(b) 
compact if it can be represented as 
$Y=\pro Y_\b$ 
where $(Y_\b)$ is a filtering projective system over $\Fschb$ 
such that each map 
$\pi^\a_{\b\b'}\,:\,Y^\a_{\b'}\to Y^\a_\b$ is affine, 

(c) 
locally compact if it can be represented as 
$Y={``\ind"}_{\a}\pro_{\b}Y^\a_\b$
where $(Y^\a_\b)$ is a bi-filtering ind-pro-system over $\Fschb$ 
such that the following conditions hold : 
for each $\b$ and $\a\leq\a'$ the structure map 
$i^{\a\a'}_\b\,:\,Y_\b^{\a}\to Y_\b^{\a'}$ is a closed embedding, 
each $\pi^\a_{\b\b'}\,:\,Y^\a_{\b'}\to Y^\a_\b$ is affine, 
for each $\a,\a'$ and $\b\leq\b'$ we have a Cartesian
commutative square 
$$\matrix 
Y^\a_{\b'}&\hookrightarrow&Y^{\a'}_{\b'}\cr 
\downarrow&&\downarrow\cr 
Y^\a_\b&\hookrightarrow&Y^{\a'}_\b.
\endmatrix$$ 
\endproclaim 

If the ind-scheme $Y$ is locally compact the ring
$\Oc(Y)=\Oc_Y(Y)$ is equal to $\pro\,\ind\,\Oc(Y^{\a}_{\b}).$
Further $Y$ can also be represented as the projective limit
$\pro Y_\beta$, where $Y_\beta={``\ind"} Y^\alpha_\beta.$
See \cite{KV1}.
For a future use set 
$Y^\alpha=\pro Y^\alpha_\beta$,
let $\pi^\a_\b: Y^\a\to Y^\a_\b$,
$\pi_\b: Y\to Y_\b$ be the projections, 
and $i^\a:Y^\a\to Y$ be the inclusion.

\vskip3mm 

\noindent{\bf (1.1.2) Examples.} 
(a) Let $V$ be a Tate vector space over $\CC$, i.e.,
a locally linearly compact linearly topological vector space. 
Then $V$ can be represented as 
$$V=\ind_\a\pro_\b V^\a_\b,$$
where $V^\a_\b$ is a bifiltering ind-pro-system of finite-dimensional 
(discrete) $k$-vector spaces such that 
all $i^{\a\a'}_\b\,:\,V^\a_\b\to V^{\a'}_\b$ are injective, 
all $\pi^\a_{\b\b'}\,:\,V^\a_{\b'}\to V^\a_{\b}$ are surjective, 
and all squares are Cartesian. 
Denoting $\underline V^\a_\b=\Spec S(V^\a_\b)^*$ the schematic affine space
corresponding to $V^\a_\b$, 
we associate to $V$ 
the locally compact ind-scheme 
$\underline V={``\ind"}\pro\underline V^\a_\b$
which is discrete if $V$ is discrete and compact if $V$ is linearly compact. 

(b) According to a theorem of Thomason, a compact scheme is the same
as a scheme that is quasicompact and quasiseparated in the usual
sense of Grothendieck. 

\proclaim{(1.1.3) Definition} 
An ind-scheme $Y$ is called locally locally compact, if 
it has an open covering by locally compact ind-schemes. 
\endproclaim 

\subhead{\bf (1.2) Smooth locally locally compact ind-schemes}\endsubhead

Recall that an ind-scheme 
$Y$ is called formally smooth if the corresponding functor $h_Y$ satisfies 
the infinitesimal lifting property of Grothendieck, i.e.,
for every ring $R$ and every nilpotent ideal $I\subset R$
the map $h_Y(R)\to h_Y(R/I)$ is surjective.

\proclaim{(1.2.1) Definition} 
We say that a locally compact ind-scheme $Y$ is smooth if it admits 
a representation as in (1.1.1) such that :

(a) 
the morphisms $\pi^\a_{\b\b'}$ are smooth,

(b)
there is an element 
$(\a,\b)$ such that $Y^\a_\b$ is smooth, 

(c) 
the ind-schemes $Y_\b$ are formally smooth. 
\endproclaim 

\proclaim{(1.2.2) Proposition} 
If $Y$ is locally locally compact smooth ind-scheme of
countable type,  and $U\subset Y$ is 
open, then $U$ is also smooth locally locally compact. 
\endproclaim

\noindent {\sl Proof :} 
Let $(Y^\alpha_\beta)$ be an ind-pro-system for $Y$ as in the definition. 
Let us
first introduce a class of open sub-ind-schemes in $Y$ that we call basic.
Let $b:\{\a\}\to\{\b\}$ be a monotone cofinal 
map such that for each $\alpha$
we have a Zariski open subset $U^\a_{b(\a)}$ in $Y^\a_{b(\a)}$ with
an equality of the preimages
$$(i^{\a\a'}_{b(\a')})^{-1} (U^{\a'}_{b(\a')})= 
(\pi_{b(\a) b(\a')}^\a)^{-1}(U^\a_{b(\a)})$$
for any $\alpha\leq\alpha'$.
Set $U^\a_\b=(\pi_{b(\a)\b}^\a)^{-1}(U^\a_{b(\a)})$ for each $\b\geq b(\a)$. 
The equality above implies that
$(i^{\a\a'}_{\b})^{-1} (U^{\a'}_{\b})=U^\a_{\b}$
for each $\b\ge b(\a')$.
Further $$U = {``\ind"}_\a \pro_{\b} U^\a_\b,$$
is an open sub-ind-scheme of $Y$, 
called a basic open sub-ind-scheme.
Our proposition follows from the lemma below.

\proclaim{(1.2.3) Lemma} 
Any open sub-ind-scheme of a locally compact ind-scheme
is a union of basic ones.
Further, a basic open sub-ind-scheme of a smooth locally compact ind-scheme
is locally compact. 
\endproclaim

\noindent {\sl Proof :} 
An open sub-ind-scheme $U$ in
$Y$ is the same as a  family of open subschemes
$U^\a$ of $Y^\a$ compatible under preimages of the ind-maps.
A basis of topology in the scheme $Y^\a$
being provided by preimages of the open sets in some $Y^\a_\b$, 
we have $(\pi^\a_{b(\a)})^{-1}(U^\a_{b(\a)})\subset U^\a$ for 
an appropriate $b(\alpha)$.
Arranging the familly $(b(\a))$ in
a monotone cofinal way gives the first assertion. 

Since $Y$ is smooth, the maps $\pi^\a_{\b\b'}$ are smooth and, 
in particular, open. 
Given $b$ and $(U^\a_{b(\a)})$ 
as above, we can  assume, by passing to a cofinal subset,
that the indexing sets are linearly ordered.
For $\b\geq b(\a)$ the open set $U^\a_\b\subset Y^\a_\b$ was already defined
as the preimage of  $U^\a_{b(\a)}$. 
For $\b\leq b(\a)$ set 
$U^\a_\b$ equal to the open set $\pi^\a_{\b b(\a)}(U^\a_{b(\a)})$.
The compatibility condition in the definition of a basic
set implies that $(U^\a_\b)$ forms a sub-ind-pro-system in 
$(Y^\a_\b)$ which satisfies
the Cartesian condition because $(Y^\a_\b)$ satisfies it. 
This finishes the proof. 
\qed

\vskip3mm

\noindent {\bf (1.2.4) Remark.} 
Let $S$ be a scheme of finite type.
The considerations of this subsection can be straightforwardly
generalized to the case of ind-schemes over $S$. In particular, we have 
the concept of a relatively smooth (locally) locally compact
ind-scheme $Y/S$ and (1.2.2) holds in the relative situation.

\vskip .3cm

\subhead{(1.3) Categories of ${\Cal O}$-modules and $\Cal D$-modules}\endsubhead

For any scheme $S$ let ${\bold O}_S$ denote the category of all 
quasi-coherent ${\Cal O}_S$-modules. 
Now, let $Y$ be an ind-scheme. 
Let $\Ob_Y$ be the category 
which is the limit of the projective system of categories 
$(\Ob_{Y^\a},(i^{\a\a'})^*).$
An object of $\Ob_{Y}$ is a system $\Fc$ of quasicoherent 
sheaves $\Fc^\alpha$ on $Y^\alpha$ together with isomorphisms 
$(i^{\alpha\alpha'})^*(\Fc^{\alpha'})\to\Fc^\alpha$
satisfying the obvious 
coherence conditions for each $\alpha\leq\alpha'\leq \alpha ''$. 

\proclaim{(1.3.1) Proposition} 
The category $\Ob_Y$ is identified 
with the category consisting of systems of quasicoherent sheaves 
$f^*(\Ec)=\Ec_f$ in $\Ob_S$ for each scheme $S$ and each 
morphism $f: S\to Y$, and isomorphisms $u^*(\Ec_{f'})\to\Ec_f$ 
given for any $u:S\to S'$ so that $f=f'\circ u.$ 
These isomorphisms are required to satisfy the compatibility condition for any 
decomposition $f= f''\circ u'\circ u$. 
\endproclaim

\proclaim {(1.3.2) Definition} 
(a)
A vector bundle on $Y$ is an object $\Ec$ in $\Ob_Y$ 
such that the quasi-coherent $\Oc_S$-module 
$f^*(\Ec)$ is locally free (possibly of infinite rank) 
for any $f:S\to Y.$ 

(b)
A pro-$\Oc$-module $\Ec$ on $Y$
is a system of pro-objects $\Ec_f$ in $\Ob_S$ 
for each $f:S\to Y$,
satisfying the compatibility conditions 
as in (1.3.1).
\endproclaim 

The category of pro-$\Oc$-modules on $Y$ is denoted by $\widehat\Ob_Y$.
If $Y$ is a scheme then a pro-$\Oc$-module on $Y$ is a pro-object in
$\Ob_Y$.
More generally
$\widehat\Ob_Y$ is the limit of the inductive system
of categories $(\widehat\Ob_{Y^\alpha},(i^{\a\a'})^*)$, i.e., to define
an object $\Ec$ in $\widehat\Ob_Y$ it suffices to define $\Ec_f$ when
$f$ is one of the canonical morphisms $i^\a:Y^\alpha\to Y$.
For such $f$ we write also $\Ec|_{Y^\alpha}$ for $\Ec_f$.

Let $\Ob^!_Y$ 
be the limit of the projective system of categories 
$(\Ob_{Y^\a},(i^{\a\a'})^!).$ 
For any scheme $S$, any ind-$S$-scheme $\rho:Y\to S$,
and any $\Ec=((i^\a)^!(\Ec))$ in $\Ob^!_Y$ the direct image $\rho_*(\Ec)$ is the
quasi-coherent sheaf over $S$ such that for each $U\subset S$ we have 
$\rho_*(\Ec)(U)$ equal to the filtered inductive limit of the
set of sections of $\Ec|_{\rho^{-1}(U)}$ supported on $Y^\a.$

Recall the embedding $Y\mapsto (Y_\roman{Zar},\Oc_Y)$ from
$\Ischb$ into a full subcategory of the category of topologically ringed spaces.
See \cite{KV2, sect.~2}.

\proclaim{(1.3.3) Proposition}
(a) The category $\widehat\Ob_Y$ is identified with the category of sheaves of
complete, separated, topological modules over the sheaf of topological rings
$\Oc_Y$.

(b) The category $\Ob_Y^!$ is identified with the full subcategory in
$\widehat\Ob_Y$ formed by the sheaves of modules whose topology is discrete.
\endproclaim

\noindent{\sl Proof :}
See \cite{BD2, 7.11.3}.
\qed

\vskip3mm

\noindent{\bf (1.3.4) Examples.}
(a) 
The cotangent sheaf of the ind-scheme $Y$ is the object $\Omega^1_Y$ 
in $\widehat\Ob_Y$
such that $\Omega^1_Y|_{Y^\a}$ is the filtered projective limit of the system
of sheaves $(\Omega^1_{Y^{\a'}}|_{Y^\a})$ with $\a'\ge\a$.

(b)
If $Y$ is discrete and formally smooth
we have the tangent bundle $\Theta_Y$ in $\Ob_Y$.
It is the unique vector bundle such that 
$\Ec \otimes\Theta_Y$ is the filtered inductive limit of the system of sheaves
$(\Homc(\Omega^1_{Y^\a},(i^\a)^!(\Ec)))$ for each $\Ec$ in $\Ob^!_Y,$ 
i.e., the vector bundle $\Theta_Y$ is dual 
to the pro-$\Oc$-module  $\Omega^1_Y$.
Further, if $f:Y\to Z$ is a formally smooth morphism of discrete
formally smooth ind-schemes then the differential of $f$ yields a surjective
map $\Theta_Y\to f^*(\Theta_Z)$.
See \cite{BD2, sect.~7.11-7.12} for details.

(c)
According to \cite{D, sect.~6.3.3}, a closed compact
subscheme $X$ of an ind-scheme $Y$
is called reasonable if for any closed subscheme $X'\subset Y$ 
containing $X$ the ideal of $X$
in $\Oc_{X'}$ is finitely generated.
Then, the normal sheaf of $X$ in $Y$ is set equal to
$N_{Y/X}=\ind N_{Y^\a/X}$ for $\a\ge \a'$ and $\a'$ such that
$X\subset Y^{\a'}$ and
$$N_{Y^\a/X}=\Der(\Oc_{Y^\a},\Oc_X)/\Der(\Oc_X,\Oc_X).$$ 
It is a quasi-coherent $\Oc_X$-module, being a filtered inductive limit
of a system of quasi-coherent $\Oc_X$-modules.

\vskip3mm

Given a scheme $S$ and a relatively smooth $S$-scheme $Z$ of finite type,
let $\Db_{Z/S}$ be the category 
of coherent right $\Cal D_{Z/S}$-modules on $Z$. 
For any relatively smooth locally compact ind-$S$-scheme $Y$
we have a double inductive system of categories 
$(\Db_{Y^\a_\b/S},i_{\b\bullet}^{\a\a'},\pi_{\b\b'}^{\a\bullet}).$
We define the category $\Db_{Y/S}$
of relative $\Dc$-modules on $Y$ to be 
the limit of this inductive system.

\subhead{\bf (1.4) The Sato Grassmannian} \endsubhead

Let $\Ec$ be a vector bundle on a discrete ind-scheme $Y$. 
We define a contravariant functor $\g_{\Ec}:{\Schb}\to {\Setb}$ 
taking $S$ to the set of pairs $(f,\Fc)$, 
with $f:S\to Y$ a morphism of ind-schemes 
and $\Fc$ a subbundle, i.e., a locally direct summand,
in $f^*(\Ec)$ of finite rank. 
The following is well-known.

\proclaim{(1.4.1) Lemma} 
The functor $\g_{\Ec}$ 
is represented by an ind-scheme $\Gr(\Ec)$ over $Y$. 
\endproclaim 

Assume now that $Y$ is a formally smooth discrete ind-scheme. 
Then $\Theta_Y$ is a vector bundle
and we have the Grassmannian 
$\Gr(\Theta_Y)$. 
For any coherent sheaf $\Fc$ on a scheme $Z$ let $\Quot_0(\Fc)$ 
be the open part of the Grothendieck Quot scheme parametrizing 
locally free quotients of $\Fc$ of finite rank. 

\proclaim{(1.4.2) Lemma} 
We have 
$\Gr(\Theta_Y)= 
``\ind"_\a``\ind"_{\a'\ge\a}\Quot_0(\Omega^1_{Y^{\a'}}|_{Y^\a}).$ 
\endproclaim 

\noindent{\sl Proof :} 
By definition, a morphism from a scheme $S$ to 
$\Gr(\Theta_Y)$ is the same as, first, a morphism 
$f:S\to Y^\a$ for some $\a$ and, second, 
an injection of a locally direct summand 
$\Ec\to f^*(\Theta_Y|_{Y^\a})$ 
where $\Ec$ is a vector bundle and injections are 
considered modulo isomorphisms of $\Ec$. 
We have 
$$f^*(\Theta_Y|_{Y^\a})= 
\ind_{\a'}f^*\Hc om_{Y^\a}(\Omega^1_{Y^{\a'}},\Oc_{Y^\a})= 
\ind_{\a'}\Hc om_S(f^*(\Omega^1_{Y^{\a'}}),\Oc_S),$$ 
where $\a'\ge\a$.
Because the functor $\Hc om$ commutes with inductive 
limits in the second argument and because $\Ec$ is locally free, we get 
$$\aligned
\Hom(\Ec,\ind_{\a'}\Hc om_S(f^*(\Omega^1_{Y^{\a'}}),\Oc_S))&= 
\ind_{\a'}\Hom(\Ec,\Hc om_S(f^*(\Omega^1_{Y^{\a'}}),\Oc_S))\\ 
&=\ind_{\a'}\Hom(f^*(\Omega^1_{Y^{\a'}}),\Ec^*),
\endaligned$$ 
where $\a'\ge\a$.
This identification takes an embedding of a locally free summand 
to an element of the limit represented by a surjection 
$f^*(\Omega^1_{Y^{\a'}})\to\Ec^*$ for some $\a'$. 
This follows from local freeness of $\Theta_Y$. 
But pairs consisting of $f:S\to Y^\a$ and surjections 
$f^*(\Omega^1_{Y^{\a'}})\to\Ec^*$ considered modulo 
isomorphisms of $\Ec^*$ form the image of $S$ by the functor represented 
by $\Quot_0(\Omega^1_{Y^{\a'}}|_{Y^\a}).$ 
\qed 

\vskip3mm 

Now assume that $Y$ is a locally compact smooth ind-scheme. 
The map $\pi_{\b\b'}:Y_{\b'}\to Y_\b$ being formally smooth 
the differential 
$d\pi_{\b\b'}:\Theta_{Y_{\b'}}\to\pi_{\b\b'}^*(\Theta_{Y_\b})$
is a surjective morphism of vector bundles. 
Hence we have the embedding of ind-schemes 
$\pi^*_{\b\b'}(\Gr(\Theta_{Y_\b}))\subset\Gr(\Theta_{Y_{\b'}})$
sending a subbundle $\pi_{\b\b'}^*(\Ec_\b)$ to 
$(d\pi_{\b\b'})^{-1}(\pi_{\b\b'}^*(\Ec_\b))$. 
We define the nonlinear Sato Grassmannian of $Y$ to be 
$$\Gr(Y)={``\ind"}\pi^*_\b\Gr(\Theta_{Y_\b}).\leqno(1.4.3)$$

\proclaim{(1.4.4) Lemma} 
The Grassmannian $\Gr(Y)$ is an ind-scheme over $Y$. 
It is independent of the choice of an ind-pro-system $(Y^\a_\b)$. 
\endproclaim 

\noindent{\sl Proof :} 
By definition and (1.4.2) we have 
$$\aligned
\Gr(Y)&={``\ind"}_\b\bigl(Y\times_{Y_\b} 
{``\ind"}_\a{``\ind"}_{\a'>\a}\Quot_0(\Omega^1_{Y^{\a'}_\b}|_{Y^\a_\b}) 
\bigr)\\
&={``\ind"}_\b{``\ind"}_\a\bigl(Y^\a\times_{Y^\a_\b} 
``\ind"_{\a'}\Quot_0(\Omega^1_{Y^{\a'}_\b}|_{Y^\a_\b})\bigr),
\endaligned$$ 
where $\a'\ge\a$.
As $Y^\a=\pro_{\b'} Y^\a_{\b'}$ for $\b'\ge\b$ and any $\b$, 
and since projective limits commute with fiber products, 
we can rewrite the above as 
$${``\ind"}_\b``\ind"_\a\pro_{\b'}\bigl(Y^\a_{\b'}\times_{Y^\a_\b} 
``\ind"_{\a'}\Quot_0(\Omega^1_{Y^{\a'}_\b}|_{Y^\a_\b})\bigr).$$ 
Next, since filtered inductive limits commute with fiber products it yields 
$${``\ind"}_\b{``\ind"}_\a\pro_{\b'}{``\ind"}_{\a'}
\bigl(Y^\a_{\b'}\times_{Y^\a_\b} 
\Quot_0(\Omega^1_{Y^{\a'}_\b}|_{Y^\a_\b})\bigr).$$ 
Furthermore, the projective limit in $\b'$ and the inductive 
limit in $\a'$ can be interchanged. 
This is because for any fixed $\a,\b$ the pro-ind-system of schemes 
$$Y^\a_{\b'}\times_{Y^\a_\b} 
\Quot_0(\Omega^1_{Y^{\a'}_\b}|_{Y^\a_\b}),$$ 
with $\a'\ge\a$, $\b'\ge\b$, has Cartesian squares, 
see \cite{KV1, (3.4.2)}. 
More precisely this pro-ind-system is the Cartesian product of a pro-system 
and an ind-system. 
So we can finally write $\Gr(Y)$ as 
$${``\ind"}_\a{``\ind"}_{\a'}{``\ind"}_\b\pro_{\b'}
\bigl(Y^\a_{\b'}\times_{Y^\a_\b} 
\Quot_0(\Omega^1_{Y^{\a'}_\b}|_{Y^\a_\b})\bigr),$$ 
where $\a'\ge \a$, $\b'\ge\b$.
For every $\a'$ the result of the remaining two limits (in $\b$ and $\b'$) 
depends only on the schemes $Y^\a$ and $Y^{\a'}$. 
Therefore the limits in $\a$, $\a'$ depend only on the ind-scheme $Y$. 
\qed 

\vskip3mm

Clearly, if $Y$ is a smooth locally locally compact ind-scheme 
we get a formally smooth ind-scheme $\Gr(Y)$ over $Y$ by glueing the
Sato Grassmannians of the pieces of a covering 
by locally compact ind-schemes. 

Let $S$ be a scheme, and 
$Y\to S$ be a relatively smooth locally locally compact ind-scheme.
We define the relative Sato Grassmannian $\Gr(Y/S)$ by replacing 
$\Theta_Y$ with $\Theta_{Y/S}$,
$\Omega^1_{Y^{\a'}}$ with $\Omega^1_{Y^{\a'}/S}$ in (1.4.1)-(1.4.2), and
$\Theta_{Y_\b}$ with $\Theta_{Y_\b/S}$ in (1.4.3). 
The relative version of (1.4.4)
is proved in the same way, yielding 
a relatively formally smooth ind-scheme
$\Gr(Y/S)$ over $Y$.

\proclaim{(1.4.5) Proposition}
(a) (base-change)
Let $Y\to S$ be a relatively smooth locally locally compact ind-scheme,
and $f:S'\to S$ be a morphism of schemes.
Then $f^*(Y)=Y\times_SS'\to S'$
is also a relatively smooth locally locally compact ind-scheme.
Further the ind-schemes $\Gr(f^*(Y)/S')$ and 
$f^*(\Gr(Y/S))=\Gr(Y/S)\times_SS'$
over $f^*(Y)$ are isomorphic.

(b) (functoriality)
Let $f:Y'\to Y$ be an isomorphism of
relatively smooth locally locally compact ind-schemes over $S$.
Then the ind-schemes $\Gr(Y'/S)$ and 
$f^*(\Gr(Y/S))=\Gr(Y/S)\times_YY'$
over $Y'$ are isomorphic.

(c) (multiplicativity)
For any relatively smooth locally locally compact ind-schemes
$Y\to S$ and $Y'\to S'$
there is a natural embedding of ind-schemes
$\Gr(Y/S)\times\Gr(Y'/S')\subset\Gr(Y\times Y'/S\times S')$
over $Y\times Y'$.
\endproclaim

\noindent{\sl Proof :}
To prove multiplicativity, notice that the ind-scheme
$Y\times Y'$ can be represented as $\pro_\b(Y_\b\times Y'_\b)$
and that 
$\Theta_{Y_\b\times Y'_\b/S\times S'}$
is equal to $\Theta_{Y_\b/S}\times\Theta_{Y'_\b/S'}$.
The rest is left to the reader.
\qed

\subhead{(1.5) The co-Sato Grassmannian}\endsubhead

Let $Y$ be a smooth locally locally compact ind-scheme. 
Then $Y$ is a reasonable formally smooth $\aleph_0$-ind-scheme in the sense of
\cite{D}. Therefore $\Omega^1_Y$ is a Tate sheaf, see \cite{D, thm. 6.2},
i.e., for any morphism $f:\Spec(R)\to Y$ the pro-$R$-module $\Omega^1_{Y,f}$
is a Tate module. In this setting Drinfeld defines 
an ind-algebraic space $\cGr(\Omega^1_Y)\to Y$ called the co-Sato Grassmannian.
By definition, a morphism $\Spec(R)\to\cGr(\Omega^1_Y)$ consists,
first, of a morphism $f:\Spec(R)\to Y$  and, second, 
of a co-projective lattice in the Tate $R$-module $\Omega^1_{Y,f}$.

\proclaim{(1.5.1) Proposition}
The ind-algebraic space $\cGr(\Omega^1_Y)$ is isomorphic to
$\Gr(Y)$ and, in particular, is an ind-scheme.
\endproclaim

\noindent{\sl Proof :}
Without loss of generality, we can assume that
$Y$ is smooth locally compact
with all $Y^\a_\b$ affine.
Let $\Kc$ be a coprojective lattice in the
pro-$\Oc$-module
$\Omega^1_Y|_{Y^\a}=
\pro\Omega^1_{Y^{\a'}}|_{Y^\a},$
for $\a'\ge\a.$
Any lattice being open in the projective limit topology,
the lattice $\Kc$ is a pull-back of some submodule
$\Kc^{\a'}\subset\Omega^1_{Y^{\a'}}|_{Y^\a}.$
Because $\Kc$ is coprojective, the quotient
$$\Fc=(\Omega^1_{Y^{\a'}}|_{Y^\a})/\Kc^{\a'}=
(\Omega^1_{Y}|_{Y^\a})/\Kc$$
is projective of finite rank over the coordinate ring of $Y^\a$.
So the choice of $\a'$ and $\Fc$ gives a morphism
$Y^\a\to\Quot_0(\Omega^1_{Y^{\a'}}|_{Y^\a}),$
and our statement follows from (1.4.2).
\qed

\subhead{(1.6) Complements on gerbes}\endsubhead

Recall that a torsor $T$ over an abelian group $A$
is a principal homogeneous $A$-space, possibly empty.
If $\phi: A\to B$ is a homomorphism of abelian groups and $T$ is an $A$-torsor,
then we have the induced $B$-torsor
$\phi_*(T) = T\otimes_A B$.

Recall that if $(\Xc, \Oc_{\Xc})$ is a ringed Grothendieck site,
then (sheaves of) locally nonempty
$\Oc_{\Xc}^\times$-torsors are in bijection with invertible $\Oc_\Xc$-modules.
If $(\Xc, \Oc_{\Xc})$ is a scheme, the latter are the same as line bundles
over $\Xc$. If $L$ is such a bundle, we denote by $L^\circ$
the corresponding ``punctured bundle'' obtained by removing the zero section.
Sections of $L^\circ$ form an $\Oc_{\Xc}^\times$-torsor. Further, if
$Y$ is an ind-scheme and $(\Xc, \Oc_{\Xc}) = (Y_{\Zar}, \Oc_Y)$, then a line
bundle $L$ on $Y$ (defined as in (1.3)) gives an $\Oc_Y^\times$-torsor still
denoted $L^\circ$.

We follow the same conventions on gerbes as in \cite{KV2, sect.~ 1}. 
This is, if $(\Xc, \Oc_\Xc)$ is a ringed site, we
distinguish two types of objects : 
$\Oc_\Xc^\times$-gerbes and sheaves of
$\Oc_\Xc^\times$-groupoids. 
If $\Cc$ is a locally connected sheaf of $\Oc_\Xc^\times$-groupoids,
then $\Cc^{\sim}$, the associated stack, 
is an $\Oc_\Xc^\times$-gerbe. 

Let $f: (\Xc, \Oc_\Xc)\to (\Yc, \Oc_\Yc)$ be a morphism of ringed sites.
It consists of 
a morphism of sites $f_\sharp: \Xc\to\Yc$,
and a morphism of sheaves of rings 
$f^\flat: f^{-1}_\sharp(\Oc_\Yc)\to\Oc_\Xc$
on $\Xc$.
If $\Lambda$ is a sheaf of $\Oc_\Yc^\times$-torsors, then we have a sheaf
of $\Oc_\Xc^\times$-torsors 
$$f^*(\Lambda)=
f_\sharp^{-1}(\Lambda)\otimes_{f_\sharp^{-1}(\Oc_\Yc^\times)}\Oc_\Xc^\times.
$$
Similarly, if $\Gc$ is an $\Oc_\Yc^\times$-gerbe (resp.\ a sheaf of
$\Oc_\Yc^\times$-groupoids), then we have an $\Oc_\Xc^\times$-gerbe
(resp.\ a sheaf of $\Oc_\Xc^\times$-groupoids) $f^*(\Gc)$.

If $\Xc$ is a ringed site and $\Gc, \Hc$ are two $\Oc_\Xc^\times$-gerbes, 
then their tensor
product $\Gc\otimes\Hc$ is the $\Oc_\Xc^\times$-gerbe such that 
$$\aligned
\Obj \bigl( (\Gc\otimes \Hc)(U)\bigr) &= 
\bigl(\Obj \, \Gc(U)\bigr) \times\bigl(\Obj\,\Hc(U)\bigr),\\
\Hom_{\Gc\otimes\Hc} \bigl( (x,y), (x', y')\bigr) &= \Hom_{\Gc}(x, x')\otimes
\Hom_\Hc(y,y').
\endaligned$$
If $\Xc_1$, $\Xc_2$ are ringed sites and 
$\Gc_i$ is an $\Oc_{\Xc_i}^\times$-gerbe, then
we have the $\Oc_{\Xc_1\times\Xc_2}^\times$-gerbe
$\Gc_1\boxtimes\Gc_2 = p_1^*(\Gc_1)\otimes p_2^*(\Gc_2).$
Here $p_i: \Xc_1\times\Xc_2\to \Xc_i$ is the obvious projection. 

If $Y$ is an ind-scheme, by an $\Oc_Y^\times$-gerbe on $Y$
we mean an $\Oc_\Yc^\times$-gerbe over the locally ringed site
$\Yc=(Y_\roman{Zar},\Oc_Y).$ 
As in (1.3.1) it is also the same as a system of $\Oc^\times_S$-gerbes
(over $S_\Zar$)
for each scheme $S$ and each morphism $f:S\to Y$ satisfying the
obvious compatibility conditions.
Taking instead $\Oc^\times_S$-gerbes over the small \'etale sites $S_\Et$
we get an \'etale $\Oc^\times_Y$-gerbe over $Y$.

\vskip3mm

\noindent {\bf (1.6.1) Groupoid ind-schemes.}
Here is a source of sheaves of $\Oc_\Xc^\times$-groupoids.
Let $Y$ be an ind-scheme. An $\Oc_Y^\times$-groupoid ind-scheme
is a system consisting of an ind-scheme $Z\to Y$ 
and a multiplicative line bundle $L$ on
$Z\times_Y Z$, i.e., a line bundle $L$ together with an isomorphism of
line bundles on $Z\times_Y Z\times_Y Z$
$$\mu: p_{12}^*(L)\otimes p_{23}^*(L)\to p_{13}^*(L),$$
satisfying the associativity condition 
and the isomorphisms of unit and inversion
$\epsilon: \Oc_Z\to \Delta^*(L),$ $i: \sigma^*(L)\to L^{-1},$
where $\Delta: Z\to Z\times_Y Z$ is the diagonal embedding and 
$\sigma: Z\times_YZ \to Z\times_Y Z$
is the permutation. 
These are required to satisfy the usual unit and inversion identities.
In particular, $i^2=\Id$.

Given an $\Oc_Y^\times$-groupoid ind-scheme $(Z\to Y,L)$, 
local sections $z: Y\to Z$ form the sheaf
of objects of a sheaf $\Zc$ of $\Oc_Y^\times$-groupoids.  
Namely, if $U$ is a scheme,
$y:U\to Y$ is a morphism, and $z_1, z_2:U\to Z$ are two sections over $y$, 
then
$\Hom_\Zc(z_1, z_2)$ is set equal to $\Gamma(U, (z_1\times z_2)^*(L^\circ)).$

We'll say that the groupoid ind-scheme $(Z\to Y,L)$ is locally trivial 
if the morphism $Z\to Y$ is locally trivial for the Zariski topology.
Then, in particular, the map $Z\to Y$ admits a section locally
on $Y_\Zar$.
Thus the sheaf of groupoids
$\Zc$ is locally non-empty and locally connected.
Therefore $\Zc^\sim$ is an $\Oc_Y^\times$-gerbe. 

If $f: Y'\to Y$ is a morphism of ind-schemes and 
$(Z\to Y,L)$ is an $\Oc_Y^\times$-groupoid ind-scheme, 
then the ind-scheme $f^*(Z) = Z\times_Y Y'\to Y'$ and the line bundle
$f^*(L)$ over $f^*(Z)\times_{Y'}f^*(Z)=f^*(Z\times_YZ)$
form an $\Oc_{Y'}^\times$-groupoid ind-scheme.
The $\Oc_{Y'}^\times$-gerbe associated to 
the groupoid ind-scheme 
$f^*(Z,L)=(f^*(Z),f^*(L))$
is equivalent to $f^*(\Zc^\sim)$. 

\vskip3mm

\noindent {\bf (1.6.2) Gerbes with connection.}
Let $\Fsetb$ be the category of nonempty finite sets and arbitrary maps. 
For every smooth algebraic variety $S$ and
every $I\in\Fsetb$ we denote by $S^{[I]}$ the formal neighborhood
of the diagonal in $S^I$. 
If $p: J\to I$ is any morphism in $\Fsetb$, we have the corresponding morphism
$\Delta_p: S^{[I]}\to S^{[J]}$,
$(x_i)\mapsto(x_{p(j)})$.
These morphisms satisfy the identity
$\Delta_{pq}=\Delta_q\circ\Delta_p$ 
for any composable pair 
$K\buildrel q\over\lra J\buildrel p\over\lra I$.

Let $Y\to S$ be an ind-scheme over $S$. 
Following Grothendieck, by an integrable connection 
on $Y$ along $S$ we mean a system 
of ind-schemes $Y_{S^{[I]}}\to S^{[I]}$
and a system of isomorphisms of ind-$S^{[I]}$-schemes 
$\nabla_{Y,p}:\Delta_p^*(Y_{S^{[J]}})\to Y_{S^{[I]}}$
such that
$\nabla_{Y, pq}=\nabla_{Y, p}\circ\Delta_p^*(\nabla_{Y,q})$
for any composable pair $p,q$.

Further, let $\Gc$ be an $\Oc^\times$-gerbe on $Y$. 
Following \cite{BM, def.~ 4.1}, by an integrable connection 
on $\Gc$ along $S$ we mean a system 
of gerbes $\Gc_I$ on $Y_{S^{[I]}}$ with $\Gc_{\{1\}}=\Gc$,
a system $\nabla_\Gc=(\nabla_{\Gc,p})$
of equivalences of gerbes 
$(\nabla_{Y,p})_*\Delta_p^*(\Gc_J)\to\Gc_I,$
and a system $\nabla_\Gc^{(2)}=(\nabla_{\Gc,p,q}^{(2)})$
of isomorphisms of equivalences
$$
\nabla_{\Gc,p}\circ(\nabla_{Y,p})_*\Delta_p^*(\nabla_{\Gc,q})\Rightarrow
\nabla_{\Gc, pq}$$
which satisfy the obvious compatibility condition for any composable triple
$p,q,r$.

\proclaim{(1.6.3) Proposition} 
Given a locally trivial $\Oc_Y^\times$-groupoid ind-scheme $(Z\to Y,L)$, 
an integrable connection on the $\Oc^\times_Y$-gerbe $\Zc^\sim$
is a pair consisting of
integrable connections on $Z$, $L$ 
which is compatible with the multiplicative structure on $L$.
The connections on $Y,Z,L$ must also be compatible in the obvious way. 
\endproclaim

For every global object $x$ of $\Gc$,
a connection
datum on $x$ is a system of objects $x_I$ of $\Gc_I$ with $x_{\{1\}}=x$
and isomorphisms
$\nabla_{\Gc,p}(x_J)\to x_I$
satisfying the compatibility conditions for any composable pair $p,q$. 

\subhead{(1.7) The determinantal gerbe} \endsubhead

Let $Y$ be a smooth locally locally compact ind-scheme and
$\Gr(Y)$ be the Sato Grassmannian. We  define the
determinantal bundle $L_Y$ on $\Gr(Y)\times_Y \Gr(Y)$ as follows.
Let $S$ be a scheme and $\phi: S\to \Gr(Y)\times_Y \Gr(Y)$ be a morphism.
So $\phi$ consists of a morphism $f: S\to Y$ and
of two collections $\Ec=(\Ec_\b)$, $\Fc=(\Fc_\b)$ of subbundles of finite rank
of $f^*(\Theta_{Y_\b})$, for large $\b$, such that 
$\Ec_{\b'}= (f^*(d\pi_{\beta\beta'}))^{-1}(\Ec_\beta)$ for $\b'> \b$ and 
similarly for
$\Fc_{\b'}$. We set 
$$(\Ec_{\b'}|\Fc_{\b'})=\Wedge^{\max}(\Ec_{\b'})
\otimes\Wedge^{\max}(\Fc_{\b'})^{-1}.$$
We have  
$(\Ec_{\b'}|\Fc_{\b'})=(\Ec_{\b''}|\Fc_{\b''})$ for each $\b'\leq\b'',$ 
yielding a line bundle $(\Ec|\Fc)$ on $S$. 
These line bundles are compatible for different $\phi$ and so yield a line
bundle $L_Y$ on $\Gr(Y)\times_Y \Gr(Y)$
such that $\phi^*(L_Y)=(\Ec|\Fc)$. 

\proclaim{(1.7.1) Proposition} 
The bundle $L_Y$ possesses a natural multiplicative structure.
Assume further that the cotangent sheaf $\Omega^1_Y$ 
is locally trivial over $Y_\Zar$.
Then $(\Gr(Y)\to Y,L_Y)$ is a locally trivial $\Oc_Y^\times$-groupoid 
ind-scheme.
\endproclaim

We'll write $\Detc_Y$ for the resulting $\Oc_Y^\times$-gerbe.

\vskip3mm

\noindent {\bf (1.7.2) Remarks.} 
(a)
According to \cite{D, sect.~3.6, 6.3.6} any smooth locally locally 
compact ind-scheme $Y$ is equipped with a determinantal 
\'etale $\Oc^\times_Y$-gerbe.
If the cotangent sheaf $\Omega^1_Y$ is locally trivial over $Y_\Zar$
then Drinfeld's determinantal gerbe is a gerbe over $Y_\Zar$,
and it is taken to the gerbe $\Detc_Y$ above
under the identification of (1.5.1).

(b)
Let $S$ be a scheme, and $Y\to S$ be a relatively
smooth locally locally compact ind-$S$-scheme.
Then the above considerations generalize to give
a construction of the determinantal line bundle $L_{Y/S}$ on
$\Gr(Y/S)\times_Y\Gr(Y/S)$ which makes $\Gr(Y/S)$ 
into an $\Oc_Y^\times$-groupoid
ind-scheme. The corresponding $\Oc_Y^\times$-gerbe will be
denoted $\Detc_{Y/S}$. 

\proclaim {(1.7.3) Proposition}
Let $Y\to S$ be a relatively
smooth locally locally compact ind-$S$-scheme
with an integrable connection $\nabla_Y$ along $S$.
Then the $\Oc_Y^\times$-gerbe 
$\Detc_{Y/S}$ inherits an integrable connection
along $S$.
\endproclaim 

\noindent{\sl Proof :}
We'll describe in details the integrable connection
on $\Gr(Y/S)$ along $S$.
The rest follows from (1.6.3) and 
is left to the reader.
We must construct a system 
of ind-$S^{[I]}$-schemes $\Gr(Y/S)_{S^{[I]}}$ 
and a system of isomorphisms of $S^{[I]}$-schemes
$\nabla_{\Gr(Y/S),p}:\Delta_p^*(\Gr(Y/S)_{S^{[J]}})\to \Gr(Y/S)_{S^{[I]}}$
satisfying the compatibility conditions above, such that
$\nabla_{\Gr(Y/S),pt}$ is equal to $\Id_{\Gr(Y/S)}$.
For each $i\in I$, to the map $p_i:pt\to I$, $pt\mapsto i$
corresponds the $i$-th projection $\Delta_i=\Delta_{p_i}:S^{[I]}\to S$.
As the ind-$S^{[I]}$-scheme $Y_{S^{[I]}}$ is isomorphic to
$\Delta_i^*(Y)$ via $\nabla_{p_i}$, it is relatively 
smooth locally locally compact by (1.4.5).
Thus the relative Sato Grassmannian $\Gr(Y_{S^{[I]}}/S^{[I]})$
is well-defined,
and we may set $\Gr(Y/S)_{S^{[I]}}=\Gr(Y_{S^{[I]}}/S^{[I]})$.
Further (1.4.5) yields an isomorphism
$$\nabla_{\Gr(Y/S),p_i}:
\Delta_i^*(\Gr(Y/S))\to
\Gr(\Delta^*_i(Y)/S^{[I]})\to
\Gr(Y/S)_{S^{[I]}}.$$
It is clear that all the required compatibility conditions are 
satisfied.
\qed

\subhead{(1.8) The sheaf of sections of a $\Cal D$-module}\endsubhead

Let $S$ be a scheme,
and $Y\to S$ be a relatively smooth locally compact ind-$S$-scheme. 
A section $\Ec$ of $\Gr(Y/S)\to Y$ gives rise to 
a functor $\Gamma_{\Ec/S}\,:\,\Db_{Y/S}\to\Ob_Y^!$ as follows. 
Recall that $\Ec$ is a family $(\Ec_\b)$ as in (1.7). 
Given $\Mc_\b^\a$ in $\Db_{Y^\a_\b/S}$, set
$\Mc^{\a'}_{\b'}$ equal to 
$\pi_{\b\b'}^{\a'\bullet}i_{\b\bullet}^{\a\a'}(\Mc^\a_\b)$
for any $\a'\geq\a$ and $\b'\geq\b$. 
By definition of the functors 
$\pi^\bullet$ on $\Cal D$-modules, we have 
$$\Mc^{\alpha'}_{\beta'} = 
\pi^{\alpha'*}_{\beta\beta'} 
i^{\alpha\alpha'}_{\beta\bullet}(\Mc^\alpha_\beta)
\otimes_{\Oc_{Y^{\alpha'}_{\beta'}}}\omega_{Y^{\alpha'}_{\beta'}/ 
Y^{\alpha'}_\beta}.$$
See \cite{KV1}.
Now, we have a similar formula : 
$$\Wedge^{\max}(\Ec_{\beta'})|_{Y^{\alpha'}_{\beta'}}=
\pi_{\beta\beta'}^{\alpha' *} 
\Bigl( \Wedge^{\max}(\Ec_\beta)|_{Y^{\alpha'}_\beta}\Bigr) \otimes 
_{\Oc_{Y^{\alpha'}_{\beta'}}}\omega^{-1}_{Y^{\alpha'}_{\beta'}/ 
Y^{\alpha'}_\beta}.$$ 
Thus for any $\a''\geq\a'\geq\a$ and $\b''\geq\b'\geq\b$ we have
$$\aligned
\Mc^{\a''}_{\b''}\otimes_{\Oc_{Y^{\a''}_{\b''}}}
\Wedge^{\max}(\Ec_{\b''})|_{Y^{\a''}_{\b''}}
&=\pi_{\b'\b''}^{\a'' *}\Bigl( 
i^{\a'\a''}_{\b'\bullet}(\Mc^{\a'}_{\b'})
\otimes_{\Oc_{Y^{\a''}_{\b'}}}
\Wedge^{\max}(\Ec_{\b'})|_{Y^{\a''}_{\b'}}\Bigr)
\\
&\supset\pi_{\b'\b''}^{\a'' *}i^{\a'\a''}_{\b' *}
\Bigl(\Mc^{\a'}_{\b'}\otimes_{\Oc_{Y^{\a'}_{\b'}}}
\Wedge^{\max}(\Ec_{\b'})|_{Y^{\a'}_{\b'}}\Bigr).
\endaligned$$ 
Therefore the vector spaces 
$$\Gamma\Bigl(Y^{\a'}_{\b'},\Mc^{\a'}_{\b'}\otimes_{\Oc_{Y^{\a'}_{\b'}}} 
\Wedge^{\max}(\Ec_{\b'})|_{Y^{\a'}_{\b'}}\Bigr)$$ 
form an inductive system. See also [KV1, (3.4.6)]. 
So, if the object $\Mc$ of $\Db_{Y/S}$ is represented by $\Mc_\b^\a$ we set 
$$\Gamma_{\Ec/S}(\Mc)=\ind_{\a'\geq\a,\b'\geq\b} 
\Gamma\Bigl(Y^{\a'}_{\b'},\Mc^{\a'}_{\b'}\otimes_{\Oc_{Y^{\a'}_{\b'}}}
\Wedge^{\max}(\Ec_{\b'})|_{Y^{\a'}_{\b'}}\Bigr).$$ 
  
\proclaim{(1.8.1) Lemma} 
The vector space 
$\Gamma_{\Ec/S}({\Cal M})$ is a discrete $\Oc(Y)$-module.
It depends only on $Y$ as an 
ind-$S$-scheme and on the objects $\Ec$ and $\Mc$ 
but not on the choice of a system $(Y^\alpha_\beta)$. 
\endproclaim 

\noindent {\sl Proof :} 
Recall that for any $\a'\geq\a$ and $\b'\geq\b$ we have
$$
\Mc^{\a'}_{\b'}\otimes_{\Oc_{Y^{\a'}_{\b'}}}
\Wedge^{\max}(\Ec_{\b'})|_{Y^{\a'}_{\b'}}
=\pi_{\b\b'}^{\a' *}
\Bigl(\Mc^{\a'}_{\b}\otimes_{\Oc_{Y^{\a'}_{\b}}}
\Wedge^{\max}(\Ec_{\b})|_{Y^{\a'}_{\b}}\Bigr).$$
Thus, for each $\alpha'$, the limit over $\beta'$
depends only on the 
scheme $Y^{\alpha'}$ and the object 
$\pi_\beta^{\alpha' *}(\Mc^{\alpha'}_\beta\otimes \Wedge^{\max} 
(\Ec_\beta)|_{Y^{\alpha'}_\beta}).$
The proof of independance finishes as in 
\cite{KV1, (4.3.3)}. 

Next, one sees that $\Gamma_{\Ec/S}(\Mc)$ possesses a natural action of
the ring $\Oc(Y).$
An element $f\in\Oc(Y)$ is a system consisting,
for each $\a'$, of a choice of a $\b'=\b'(\a')$ and of a function 
$f^{\a'}_{\b'}\in\Oc(Y^{\a'}_{\b'})$ 
satisfying natural compatibility conditions.
Here we can assume that $\a'\ge\a$ and $\b'\ge\b$.
To spell out these conditions, we can assume that $\b'$ depends on $\a'$
in a monotone way. 
Then, taking any $\a\le\a'\le\a''$ and the corresponding $\b\le\b'\le\b''$,
the condition writes
$(i^{\a'\a''}_{\b''})^*(f^{\a''}_{\b''})=
\pi^{\a' *}_{\b'\b''}(f^{\a'}_{\b'}).$
Recall that we have an inclusion
$$
i^{\a'\a''}_{\b'' *}
\pi_{\b'\b''}^{\a' *}
\Bigl(\Mc^{\a'}_{\b'}\otimes_{\Oc_{Y^{\a'}_{\b'}}}
\Wedge^{\max}(\Ec_{\b'})|_{Y^{\a'}_{\b'}}\Bigr)
\subset
\Mc^{\a''}_{\b''}\otimes_{\Oc_{Y^{\a''}_{\b''}}}
\Wedge^{\max}(\Ec_{\b''})|_{Y^{\a''}_{\b''}}
.$$
Thus the action of
$f^{\a'}_{\b'}$ on
$\Gamma\bigl(Y^{\a'}_{\b'},
\Mc^{\a'}_{\b'}\otimes\Wedge^{\max}(\Ec_{\beta'})
|_{Y^{\alpha'}_{\beta'}}\bigr)$
is compatible with the action of
$f^{\a''}_{\b''}$ on
$\Gamma\bigl(Y^{\a''}_{\b''},
\Mc^{\a''}_{\b''}\otimes\Wedge^{\max}(\Ec_{\b''})|_{Y^{\a''}_{\b''}}\bigr)$
under the natural map of the first space to the second one.
So $\Gamma_{\Ec/S}(\Mc)$ is an $\Oc(Y)$-module.
It is discrete since it is represented as an inductive limit 
of discrete vector spaces.
\qed

\vskip3mm

Let now $Y\to S$ be a relatively
smooth locally locally compact ind-$S$-scheme,
and $\Mc$ be as before. 
Localizing on $U\subset Y$ we get a functor 
$\Gamma_{(\Ec|_U)/S}\,:\,\Db_{U/S}\to\Ob_U^!$.

\proclaim{(1.8.2) Proposition}
For any object $\Ec$ of $\Detc_{Y/S}(Y)$
the relative global sections yield a functor  
$\Gamma_{\Ec/S}:\Db_{Y/S}\to\Ob_Y^!$.
Any isomorphism $\Ec\to\Ec'$ in $\Detc_{Y/S}(Y)$
yields an isomorphism of functors
$\Gamma_{\Ec/S}\to\Gamma_{\Ec'/S}.$
For every $f\in\Oc(Y)^\times$ the corresponding
automorphism of $\Ec$ is taken to the multiplication by $f$ on 
$\Gamma_{\Ec/S}.$
\endproclaim

\noindent{\sl Proof :}
First, observe that for any locally compact open ind-$S$-schemes $V\subset U$ 
we have a morphism of functors
$\Gamma_{(\Ec|_U)/S}\to\Gamma_{(\Ec|_V)/S}$ by (1.8.1).
Fix an open covering $Y=\bigcup_iU_i$ by locally compact ind-$S$-schemes.
For each $i$ we have the functor
$\Gamma_{(\Ec|_{U_i})/S}\,:\,\Db_{U_i/S}\to\Ob_{U_i}^!.$
For each $i,j$ the open ind-subscheme $U_i\cap U_j$ 
is  locally locally compact by (1.2.2).
Fix an open covering $U_i\cap U_j=\bigcup_k U_{ij,k}$ 
by locally compact ind-$S$-schemes.
For each object $\Mc\in\Db_{Y/S}$ we set
$$\Gamma_{\Ec/S}(\Mc)=\roman{\Ker}\Bigl\{
\prod_i\Gamma_{(\Ec|_{U_i})/S}(\Mc|_{U_i})\to
\prod_{i,j,k}\Gamma_{(\Ec|_{U_{ij,k}})/S}(\Mc|_{U_{ij,k}})\Bigr\}.
$$
The action of the functor $\Gamma_{\Ec/S}$ on morphisms is defined similarly.
\qed

\vskip3mm

Given two relatively smooth locally locally compact ind-schemes
$Y\to S$ and $Y'\to S'$,
let $\boxtimes$ denote the external tensor products
$$\Db_{Y/S}\times\Db_{Y'/S'}\to\Db_{Y\times Y'/S\times S'},
\quad
\Ob^!_{Y}\times\Ob^!_{Y'}\to\Ob^!_{Y\times Y'}.$$
The compatibility of the functor $\Gamma$ with Cartesian products is expressed
by the following proposition.

\proclaim{(1.8.3) Proposition}
The assignement $(\Ec,\Ec')\mapsto\Ec\times\Ec'$
yields an equivalence of $\Oc^\times_{Y\times Y'}$-gerbes
$\Detc_{Y/S}\boxtimes\Detc_{Y'/S'}\to\Detc_{Y\times Y'/S\times S'}$.
For any $\Mc$, $\Mc'$ the space of sections 
$\Gamma_{\Ec\times\Ec'/S\times S'}(\Mc\boxtimes\Mc')$
is equal to 
$\Gamma_{\Ec/S}(\Mc)\boxtimes\Gamma_{\Ec'/S'}(\Mc').$
\endproclaim

\noindent{\sl Proof :}
Let $i_{Y,Y'}$ be the embedding of
$\Gr(Y/S)\times\Gr(Y'/S')$ into $\Gr(Y\times Y'/S\times S')$.
We have an isomorphism 
$L_{Y/S}\boxtimes L_{Y'/S'}=
(i_{Y,Y'}\times i_{Y,Y'})^*(L_{Y\times Y'/S\times S'}),$
yielding an equivalence of $\Oc^\times_{Y\times Y'}$-gerbes
$\Detc_{Y/S}\boxtimes\Detc_{Y'/S'}\to\Detc_{Y\times Y'/S\times S'}$.
\qed

\vskip3mm

Let $Y\to S$ be a relatively smooth locally locally compact ind-scheme,
and $f:S'\to S$ be a morphism of schemes.
Compatibility of the functor $\Gamma$ with base change
is expressed by the following proposition.

\proclaim{(1.8.4) Proposition}
The assignement $\Ec\mapsto f^*(\Ec)$ yields an equivalence of
$\Oc^\times_{f^*(Y)}$-gerbes $f^*(\Detc_{Y/S})\to\Detc_{f^*(Y)/S'}$.
For each $\Mc$ we have
$\Gamma_{f^*(\Ec)/S'}(f^*(\Mc))=f^*(\Gamma_{\Ec/S}(\Mc)).$
\endproclaim

\noindent{\sl Proof :}
There is an isomorphism
$$f^*(\Gr(Y/S))=
``\ind"\pi_\b^*(\Gr(\Theta_{f^*(Y_\b)/S'}))=
\Gr(f^*(Y)/S').$$
So the groupoid ind-schemes
$f^*\bigl(\Gr(Y/S),L_{Y/S}\bigr)$
and
$\bigl(\Gr(f^*(Y)/S'),L_{f^*(Y)/S'}\bigr)$
are isomorphic.
The rest is left to the reader.
\qed

\vskip3mm

Functoriality of $\Gamma$, in the sense of (1.4.5),
is left to the reader.

\vfill\eject

\centerline{\bf 2. Factorization structures.}

\vskip 1cm

\subhead {(2.1) Coherence on Cartesian powers}\endsubhead

The concept of a factorization algebra as defined in [BD1]
and used in \cite{KV1} is based on a combinatorial formalism involving
objects defined on various Cartesian powers $C^I$ for a given curve $C$.
In this paper we need a further level of coherence for this
formalism since we work with gerbes. So we start with
an equivalent reformulation of the formalism better adapted for
studying coherence.

Let $\Fsetb_0$ be the category of nonempty finite sets and their surjections. 
It has a final object $\{1\}$ (a one-point set)
and a monoidal structure $\coprod$ (disjoint union)
but no unit object for $\coprod$. 
If $p: J\to I$ and $p': J'\to I'$
are two morphisms of $\Fsetb_0$, we denote by 
$p\coprod p': I\coprod I'\to J\coprod J'$
their disjoint union.

Let $C$ be a smooth algebraic curve. 
For every morphism $p: J\to I$ in $\Fsetb_0$
we denote by  $C^p$ the open subset in $C^J$ 
consisting of the $J$-uples $(c_j)$
such that $c_j\neq c_{j'}$ for $p(j)\neq p(j')$. 
We'll write $p_J$, or simply $J$, for the unique map $J\to\{1\}$.
Notice that $C^{p_J}=C^J$.
We'll also write $1_J:J\to J$ for the identity.

Let $K\buildrel q\over\lra J\buildrel p\over\lra I$ be a composable pair of
surjections. We have the diagonal map 
$\Delta_{p,q}: C^p\to C^{pq},$
$(c_j)\mapsto(c_{q(k)}),$
which is a closed embedding, and the off-diagonal map
$j_{p,q}: C^{q}\to C^{pq},$
$(c_k)\mapsto(c_k),$
which is an open embedding. 
For each $p,p'$ we have also the map 
$i_{p,p'}:C^{p\coprod p'}\to C^{p}\times C^{p'},$
$(c_k)\mapsto(c_k),$
which is an open embedding. 
The maps above fit into the following commutative diagrams,
for any composable triple 
$L\buildrel r\over\lra K\buildrel q\over \lra J\buildrel p\over
\lra I$ of surjections :
$$\matrix  C^{p}&\buildrel\Delta_{p,q}\over\lra& C^{pq}&\cr
\Delta_{p, qr}&\searrow&\big\downarrow&\Delta_{pq,r}\cr
&&C^{pqr}&\endmatrix, \leqno (2.1.1)$$

$$\matrix C^{r}& \buildrel j_{q,r}\over\lra&C^{qr}&\cr
j_{pq,r}&\searrow&\big\downarrow&j_{p, qr}\cr
&&C^{pqr}&\endmatrix, \leqno (2.1.2)$$

$$\matrix & C^{q}&\buildrel j_{p,q}\over\lra& C^{pq}&\cr
\Delta_{q,r}&\big\downarrow&&\big\downarrow&\Delta_{pq,r}\cr
&C^{qr}&\buildrel j_{p, qr}\over\lra&C^{pqr}&\endmatrix . \leqno (2.1.3)$$
Given $i\in \{ 1, ..., n\}$ and  
a composable chain of surjections 
$$I_n\buildrel p_n\over\lra I_{n-1}\buildrel p_{n-1}\over\lra ... \buildrel
p_2\over\lra I_1\buildrel p_1\over\lra I_0 \leqno (2.1.4)$$
we have a diagram formed by the varieties
$$C^{p_a p_{a+1} ... p_b}, \quad 1\leq a\leq i\leq b\leq n,\leqno(2.1.5)$$
and various instances of maps $\Delta_{p,q}$, $j_{p,q}$ among them.
This diagram is the 1-skeleton of the product of oriented simplices
$\Delta^{i-1}\times\Delta^{n-i}$. Further, any 2-face of this 
$\Delta^{i-1}\times\Delta^{n-i}$ is a diagram of the form
(2.1.1)-(2.1.3) and so is commutative. Therefore the entire diagram
is commutative. In other words, we have the ``coherence theorem'':

\proclaim{(2.1.6) Proposition} In the situation (2.1.4) the various
instances of $\Delta$ and $j$ define a unique map 
$C^{p_i}\to C^{p_1 ... p_n}$.
\endproclaim

\proclaim{(2.1.7) Definition} 
We will refer to the $\Delta^{i-1}\times\Delta^{n-i}$-diagram
in the situation of (2.1.4)-(2.1.5) as the diagram of the type
$(p_i\to p_1... p_n)$. 
\endproclaim

\subhead {(2.2) Factorization semigroups}\endsubhead

We keep the notations and conventions of (2.1).

\proclaim{(2.2.1) Definition} 
Let $Y_C\to C$ be a formally smooth ind-$C$-scheme 
equipped with an integrable connection along $C$. 
A factorization semigroup on $Y_C$ is a system consisting of :

(a) formally smooth ind-schemes $\rho_p:Y_p\to C^{p}$ 
equipped with integrable connections along $C^{p}$, 
so that $Y_{\{1\}}=Y_C$,

(b) for any composable pair $p,q$
in $\Fsetb_0$, isomorphisms of
relative ind-schemes with connections
$\varkappa_{p,q}: \Delta_{p,q}^* (Y_{pq})\to Y_p$
and
$\kappa_{p,q}: j_{p,q}^* (Y_{pq})\to Y_q$
satisfying the compatibility conditions lifting (2.1.1)-(2.1.3) :
$$\varkappa_{p, qr} = \varkappa_{p,q}\circ\Delta_{p,q}^*(\varkappa_{pq,r}): 
\Delta_{p, qr}^* (Y_{pqr})\to Y_p,
$$
$$\kappa_{pq,r}= 
\kappa_{q,r}\circ j_{q,r}^*(\kappa_{p, qr}): j_{pq,r}^*(Y_{pqr})\to Y_r,
$$
$$
\kappa_{p,q}\circ j_{p,q}^*(\varkappa_{pq,r})=\varkappa_{q,r}\circ\Delta_{q,r}^*
(\kappa_{p, qr}): j_{p,q}^*\Delta_{pq,r}^*(Y_{pqr})=
\Delta_{q,r}^* j_{p, qr}^*(Y_{pqr}) \to Y_q,
$$

(c) 
for any pair $p,p'$ in $\Fsetb_0$, isomorphisms  
$\sigma_{p,p'}:i_{p,p'}^*(Y_p\times Y_{p'})\to Y_{p\coprod p'}$.
\endproclaim

\proclaim{(2.2.2) Definition} 
A factorization semigroup $(\rho_p:Y_p\to C^p)$ is said to be commutative
if the maps $\varkappa$, $\kappa$ factor through a morphism
of $C^{\{1,2\}}$-schemes
$Y_{\{1,2\}}\to Y_{\{1\}}\times Y_{\{2\}}$.
\endproclaim

\noindent {\bf (2.2.3) Example.} 
The collection $(C^p)$ forms a commutative factorization semigroup which
we call the unit semigroup.

\vskip3mm

%
%

\noindent {\bf (2.2.4) Remarks.} 
(a)
The definition (2.2.1) is equivalent to \cite{KV1, (2.2.1)}.
Indeed, given a system $(Y_p)$ as before,
we define $Y_I = Y_{p_I}$. 
Then the $Y_I$ satisfy the conditions of \cite{loc.\ cit.}.
Conversely, given $(Y_I)$ as in \cite{loc.\ cit.} and $p: J\to I$ a surjection,
we define $Y_p= j_{p_I,p}^*(Y_J)$.
Then the $Y_p$ satisfy the conditions of (2.2.1). 

(b)
The map $Y_{\{1,2\}}\to Y_{\{1\}}\times Y_{\{2\}}$ in (2.2.2)
is opposite to the map 
in \cite{BD1, 3.10.16}
in the axioms of commutative chiral monoids.

(c)
In 
\cite{BD1, 3.10.16}
the authors impose that the closure in
$Y_I$ of the complement to the preimage of the discriminant divisor in
$C^I$ equals $Y_I$. 
This condition is weaker than the formal smoothness
of the map $Y_I\to C^I$.
For formal loop spaces this formal smoothness is proved in \cite{KV1, (2.6.2)}.

\subhead {(2.3) Factorization semigroups of local nature}\endsubhead

We denote by $\gen$ the Lie algebra $\Der\, \CC[[t]]$ and by
$K$ the group-subscheme of the group-ind-scheme
$$\Aut\,\CC[[t]]=
\ind_n\Spec\bigl(\CC[a_0,a_1^{-1},a_1,a_2,a_3,...]/(a_0^n)\bigr)$$
defined by $a_0=0$.
So for a ring $R$ an $R$-point of $K$ is a formal change of
coordinates
$t\mapsto a_1 t + a_2 t^2 + ...$ with
$a_1\in R^\times$ and $a_i\in R$ for $i\geq 2.$
The Lie algebra $\gen$ and the group scheme $K$ form a Harish-Chandra pair.
By an action of $(\gen, K)$ on an ind-scheme $Y$ we mean an action of $K$
by automorphisms and an action of $\gen$ by derivations (infinitesimal
automorphisms) which are compatible. 

Let $C$ be as before and $\widehat {C}\to C$ be the scheme whose points are
pairs $(c, t_c)$ where $c$ is a point of $C$ and $t_c$ is a formal coordinate
near $c$. The Harish-Chandra pair $(\gen, K)$ acts on $\widehat{C}$ with the
action of $K$ preserving the projection $\widehat{C}\to C$ and the action
of the element $d/dt$ of $\gen$ defining an integrable
connection on $\widehat{C}$ along $C$.

Let $Y$ be an  ind-scheme with a $(\gen, K)$-action. 
We form the ind-scheme
$Y_C= Y\times_K \widehat{C}$
over $C$ which inherits an integrable
connection along $C$ and is formally smooth over $C$
if $Y$ is formally smooth. A structure
of a factorization semigroup on $Y_C$ will be referred to as a structure of
a factorization semigroup on $Y$. We will call such factorization semigroups
{\it factorization semigroups of local nature}.

\subhead{\bf (2.4) Factorization algebras and vertex algebras}\endsubhead

\proclaim{(2.4.1) Definition} 
Let $\Ec_C$ be a quasi-coherent sheaf over $C$
equipped with a connection $\nabla_{\Ec_C}$ along $C$. 
A factorization algebra on $\Ec_C$ is a system consisting of :

(a) 
quasi-coherent sheaves $\Ec_p$ over $C^{p}$
equipped with integrable connections $\nabla_{\Ec_p}$, 
so that $(\Ec_{\{1\}},\nabla_{\Ec_{\{1\}}})=(\Ec_C,\nabla_{\Ec_C})$,

(b) for any composable pair $p,q$
in $\Fsetb_0$, isomorphisms of
quasicoherent sheaves with connections
$\mu_{p,q}: \Delta_{p,q}^* (\Ec_{pq})\to \Ec_p$,
$\lambda_{p,q}: j_{p,q}^* (\Ec_{pq})\to \Ec_q$
satisfying the compatibility conditions lifting (2.1.1)-(2.1.3) :
$$\mu_{p, qr} = \mu_{p,q}\circ\Delta_{p,q}^*(\mu_{pq,r}): 
\Delta_{p, qr}^* (\Ec_{pqr})\to \Ec_p,
$$
$$\lambda_{pq,r} = \lambda_{q,r}\circ j_{q,r}^*(\lambda_{p, qr}): 
j_{pq,r}^*(\Ec_{pqr})\to \Ec_r,
$$
$$
\lambda_{p,q}\circ j_{p,q}^*(\mu_{pq,r})= \mu_{q,r}\circ\Delta_{q,r}^*
(\lambda_{p, qr}): j_{p,q}^*
\Delta_{pq,r}^*(\Ec_{pqr})
=\Delta_{q,r}^* j_{p, qr}^*(\Ec_{pqr}) \to \Ec_q,
$$

(c) 
for any pair $p,p'$ in $\Fsetb_0$,
isomorphisms  
$\nu_{p,p'}:
i_{p,p'}^*(\Ec_p\boxtimes \Ec_{p'})\to\Ec_{p\coprod p'}$,

(d)
a covariantly constant 
section $\zeta$ of $\Ec_C\to C$ such that for any local section $f$ 
the product $\zeta\boxtimes f$ is a local section of
$\Ec_{\{1,2\}}\to C^{\{1,2\}}$ 
such that $\Delta^*_{\{1\},\{1,2\}}(\zeta\boxtimes f)=f$.
\endproclaim

As before, this definition is equivalent to \cite{BD1, sect.~3.4}. 
Indeed, given a system $(\Ec_p)$
we define $\Ec_I = \Ec_{p_I}$. 
Then the system $(\Ec_I)$ satisfies the conditions of \cite{loc.~ cit.}.
Conversely, given $(\Ec_I)$ as in \cite{loc.~cit.} and $p: J\to I$ a surjection,
we define $\Ec_p= j_{p_I,p}^*(\Ec_J)$.
Then the system $(\Ec_p)$ satisfies the conditions of (2.4.1). 

\proclaim{(2.4.2) Definition} 
The factorization algebra $(\Ec_p)$
is commutative if the isomorphisms
$\l$, $\nu$ 
factor through a morphism
$\Ec_{\{1\}}\boxtimes\Ec_{\{2\}}\to\Ec_{\{1,2\}}$.
\endproclaim

\noindent{\bf (2.4.3) Examples.}
(a)
The system of coherent sheaves $(\Oc_{C^p})$
is a commutative factorization algebra on $\Oc_C,$
called the trivial factorization algebra.
Given $(\Ec_p)$ as above, from \cite{BD1, (3.4.4)} there is an unique morphism
of factorization algebras $(\zeta_p:\Oc_{C^p}\to\Ec_{C^p})$
such that $\zeta_{\{1\}}=\zeta$.

(b)
Let $(\rho_p:Y^0_p\to C^p)$ 
be a factorization semigroup consisting of schemes
(possibly of infinite type). 
Then the system of quasi-coherent sheaves $(\rho_{p*}(\Oc_{Y^0_p}))$
is a factorization algebra (with the units equal to the constant functions 1).
If the factorization semigroup $(Y^0_p)$ is commutative, 
so is also the factorization algebra $(\rho_{p*}(\Oc_{Y^0_p}))$.

(c)
Let $(\rho_p:Y_p\to C^p)$ be a factorization semigroup on $Y_C$ consisting of
smooth locally locally compact ind-schemes, and
$(Y^0_p)$ be a factorization sub-semigroup consisting of smooth 
(compact) reasonable subschemes.
Then the normal sheaf $N_{Y_p/Y^0_p}$ is well-defined for each $p$,
and it is a vector bundle (of infinite rank) over $Y^0_p$.
The system of symmetric algebras
$S_{Y_p/Y^0_p}=S_{\Oc_{Y^0_p}}(N_{Y_p/Y^0_p})$
yields the factorization algebra 
$(\rho_{p*}(S_{Y_p/Y^0_p}))$
which may not be commutative,
even if the factorization semigroup $(Y^0_p)$ is commutative.
We'll see in (3.2.3) that it is indeed commutative for formal loop spaces,
although, in this case, $(Y_p)$ is not a commutative semigroup.

\vskip3mm

Recall that a $\ZZ_{\ge 0}$-graded vertex algebra is a 
$\ZZ_{\ge 0}$-graded $\CC$-vector space
$V=\bigoplus_{n\ge 0}V_n$ equipped with a distinguished element $1$ in $V_0$,
an endomorphism $\partial$ of $V$ of degree one, and a family of
bilinear operations
$(x,y)\mapsto x\circ_n y$, such that 
$V_i\circ_n V_j\subset V_{i+j-n-1}$,
subject to the axioms in \cite{GMS2, def.~ 0.4}.
One defines $Y(x,t)$ to be the operator formal series
$y\mapsto\sum_n(x\circ_ny)t^{-n-1}.$
Further, a morphism of $\ZZ_{\ge 0}$-graded vertex algebras 
$T:V\to W$ is a homogeneous linear map such that 
$T(1)=1$
and
$T(Y(x,t)y)=Y(Tx,t)Ty$
for each $x,y\in V.$

\vskip3mm

\noindent{\bf (2.4.4) Example.}
 If $A$ is a commutative algebra with a derivation $\partial$, 
it is made into a vertex algebra with
$Y(a,t)b=\exp(t\,\partial)(a)\cdot b.$
Such vertex algebras will be called commutative (or holomorphic).
They are characterized by the property that
$Y(a,t)$ belongs to $\End(A)[[t]]$.
Note that $a\circ_{-1}b=a\cdot b$ is the usual product in $A$.

\vskip3mm

Let $V$ be a vector space with a $(\gen, K)$-action. The action of the
subgroup $\GG_m\subset K$ defines a $\ZZ$-grading on $V$. Let us assume
that this grading is concentrated in degrees $\ge 0$. 
As in (2.3), we associate to $V$ a quasi-coherent sheaf
$\Vc = V\times_K \widehat{C}$
on $C$ with a connection. 
One of the main results of \cite{BD1} can be formulated as follows, 
see in particular \cite{loc.~cit., (3.4.20)}.

\proclaim{(2.4.5) Proposition}
A factorization algebra structure on $\Vc$ gives rise to a structure of
a ($\ZZ_{\geq 0}$-graded) vertex algebra on $V$. 
The vertex algebra $V$ is commutative iff the factorization
algebra $\Vc$ is commutative.
\endproclaim

\subhead{(2.5) Factorizing functions}\endsubhead

Let $(Y_p)$ be a factorization semigroup on a 
formally smooth ind-$C$-scheme $Y_C$ with a connection 
$\nabla_{Y_C}$ along $C$.

\proclaim{(2.5.1) Definition}
An invertible  function $f_C$ in $\Oc^\times(Y_C)$ is called (multiplicatively)
factorizable if 
it is covariantly constant with respect to $\nabla_{Y_C}$,
and there exist invertible functions $f_p$ in $\Oc^\times(Y_p)$
satisfying the conditions 
$\Delta_{p,q}^*(f_{pq})=\varkappa_{p,q}^*(f_p),$
$j_{p,q}^*(f_{pq})=\kappa_{p,q}^*(f_q),$
$i_{p,p'}^*(f_p\times f_{p'})=\sigma_{p,p'}^*(f_{p\coprod p'})$
for any pairs $p,q$ and $p,p'$ as above.
\endproclaim

Note that, for each function $f_C$, the system $(f_p)$ is unique if it exists. 
We denote by $\Oc^\times(Y_C)_{fact}$ the multiplicative group of
factorizable functions. 

Assume that $(Y_p)$ is a factorization semigroup of local nature on $Y_C$, 
coming from a $(\gen, K)$-ind-scheme $Y$. 
Then the localization construction of (2.3) establishes an isomorphism
$\Oc^\times (Y)^\gen \to \Oc^\times(Y_C)^\nabla,$
$f\mapsto f_C,$
where the RHS is the group of covariantly constant functions with respect to
$\nabla_{Y_C}$. So we will
say that a function
$f$ in $\Oc^\times(Y)^\gen$ is factorizable if $f_C$ satisfies the
conditions of (2.5.1) and we will denote $\Oc^\times(Y)_{fact}$
the group formed by such functions. 

\vskip3mm

\noindent{\bf (2.5.2) Remark.}
Let $Y$ be a smooth locally locally compact $(\gen,K)$-ind-scheme.
Then $Y_C$ is also smooth locally locally compact.
Thus the Sato Grassmannian $\Gr(Y)$ and the relative Sato Grassmannian 
$\Gr(Y_C/C)\to C$ are well-defined.
Further $\Gr(Y)$ is equipped with a $(\gen,K)$-action such that
$\Gr(Y_C/C)=\Gr(Y)\times_K\widehat C$,
and (1.7.3) yields a connection on
$\Gr(Y_C/C)$ along $C$ inducing a bijection
$\Gamma(Y,\Gr(Y))^\gen\to\Gamma(Y_C,\Gr(Y_C/C))^\nabla.$

\subhead{(2.6) Factorizing line bundles and sections}\endsubhead

Let $(Y_p)$ be a factorization semigroup on a 
formally smooth ind-$C$-scheme $Y_C$ with a connection 
$\nabla_{Y_C}$ along $C$.
Let $L_C$ be a line bundle over $Y_C$
with a connection $\nabla_{L_C}$ along $C$ which is compatible with
$\nabla_{Y_C}$.

\proclaim{(2.6.1) Definition}
A factorizing line bundle structure on $L_C$ is a system consisting of :

(a) line bundles $(L_p,\nabla_{L_p})$ over $Y_p$ with
integrable connections along $C^{p}$ compatible with $\nabla_{Y_p}$,
so that
$(L_{\{1\}},\nabla_{L_{\{1\}}})=(L_C,\nabla_{L_C})$,

(b) for any composable pair 
$p,q$ in $\Fsetb_0$,
isomorphisms $u_{p,q}: 
\Delta_{p,q}^*(L_{pq})\to\varkappa_{p,q}^*(L_p)$
and $v_{p,q}: 
j_{p,q}^*(L_{p,q})\to\kappa_{p,q}^*(L_q)$
satisfying obvious compatibility conditions,

(c) for any pair 
$p,p'$ in $\Fsetb_0$,
isomorphisms
$w_{p,p'}:
i_{p,p'}^*(L_p\boxtimes L_{p'})\to\sigma_{p,p'}^*(L_{p\coprod p'}).$
\endproclaim 

\proclaim {(2.6.2) Definition}
Let $(L_p)$ be a factorizing line bundle on $L_C$. 
A nonvanishing section $s_C$ in $\Gamma(Y_C, L_C)$
is called factorizable, if 
it is covariantly constant with respect to $\nabla_{L_C}$
and there exist nonvanishing sections $s_p$ of $L_p$
satisfying the following conditions,
for any pairs $p,q$ and $p,p'$ as above :
$$u_{p,q}(\Delta_{p,q}^*(s_{pq}))=\varkappa_{p,q}^*(s_p),
\quad
v_{p,q}(j_{p,q}^*(s_{pq}))=\kappa_{p,q}^*(s_q),$$
$$w_{p,p'}(i^*_{p,p'}(s_p\boxtimes s_{p'}))=
\sigma_{p,p'}^*(s_{p\coprod p'}).$$
\endproclaim

As in (2.5) the system $(s_p)$ is unique if it exists.
We denote by $\Gamma(Y_C, L^\circ_C)_{fact}$ the set of factorizing sections.

\proclaim{(2.6.3) Proposition}
Assume that $L_C$ comes from a line bundle $L$ on the ind-scheme $Y$.
The set $\Gamma(Y_C, L^\circ_C)_{fact}$ is a torsor (possibly empty)
over the group $\Oc^\times(Y)_{fact}$.
\endproclaim

\noindent{\bf (2.6.4) Remark.}
Once again, the definition of factorizing line bundles is similar to the one
in \cite{BD1, (3.10.16)}.

\subhead{(2.7) Factorization gerbes}\endsubhead

Let $(Y_p\to C^{p}, \nabla_{Y_p}, \varkappa_{p,q}, \kappa_{p,q})$ 
be a factorization semigroup on $Y_C$. 

\proclaim{(2.7.1) Definition}
A factorization $\Oc^\times$-gerbe over $(Y_p)$ is a datum consisting of :

(a) for any $p$ in $\Fsetb_0$,
an $\Oc^\times_{Y_p}$-gerbe $\Gc_p$
with an integrable connection along $C^{p}$, 

(b) for any composable pair 
$p,q$ in $\Fsetb_0$,
equivalences of gerbes with connections
$A_{p,q}: \Delta^*_{p,q}(\Gc_{pq})\to\varkappa_{p,q}^*(\Gc_p)$ and
$B_{p,q}: j^*_{p,q}(\Gc_{pq})\to\kappa^*_{p,q}(\Gc_q),$

(c) for any pair $p, p'$ in $\Fsetb_0$, 
an equivalence of $\Oc^\times$-gerbes with connections
$C_{p, p'}: 
i^*_{p,p'}(\Gc_p\boxtimes\Gc_{p'})\to
\sigma_{p,p'}^*(\Gc_{p\coprod p'}),$

(d) for any composable triple 
$p,q,r$ in $\Fsetb_0$,
isomorphisms of equivalences
$$
\a_{p,q,r}: 
(\Delta_{p,q}^*(\varkappa_{pq,r}))^*(A_{p,q})\circ\Delta_{p,q}^*(A_{pq,r})
\Rightarrow A_{p, qr},
$$
$$
\b_{p,q,r}: 
(j^*_{q,r}(\kappa_{p,qr}))^*(B_{q,r})\circ j_{q,r}^*(B_{p, qr})
\Rightarrow B_{pq, r},
$$
$$
\g_{p,q,r}:
(j_{p,q}^*(\varkappa_{pq,r}))^*(B_{p,q})\circ j_{p,q}^*(A_{pq, r})\Rightarrow 
(\Delta_{q,r}^*(\kappa_{p,qr}))^*(A_{q,r})\circ\Delta_{q,r}^*(B_{p, qr}),
$$

(e) associativity and commutativity constraints for the $C_{p, p'}$.

\noindent
The data (a)-(e) are required to satisfy the following coherence conditions :

(f) for any composable  4-tuple $p,q,r,s$ in $\Fsetb_0$,
the isomorphisms $\a, \b, \g$ 
fit into four commutative polytopes : 
two 3-simplices of types $(p\to pqrs)$ and $(q\to pqrs)$
and two prisms $\Delta^1\times\Delta^2$ 
of types $(q\to pqrs)$ and $(r\to pqrs)$.
Here the commutativity of 3-dimensional polytopes in the 2-category of
gerbes is understood in the sense of pasting, see \cite{M} for
background,

(g) the standard coherence conditions for commutativity and associativity
of the equivalences $(C_{p,p'})$. 
\endproclaim

The following example of factorization gerbe is essential
for the rest of the paper.

\proclaim{\bf (2.7.2) Proposition}
Assume that each $Y_p$ is smooth locally locally compact over $C^p$
with a locally trivial cotangent sheaf.
Then we have the relative determinantal $\Oc_{Y_p}^\times$-gerbe 
$\Detc_{Y_p/C^{p}}$. 
The system $(\Detc_{Y_p/C^{p}})$ 
forms a factorization gerbe. 
\endproclaim

\noindent{\sl Proof :}
Observe that the system of ind-schemes $(\Gr(Y_p/C^{p}))$
is not a factorization semigroup, because 
there is no isomorphism of ind-schemes 
$$i^*_{p,p'}\bigl(\Gr(Y_{p}/C^{p})\times\Gr(Y_{p'}/C^{p'})\bigr)\to
\Gr(Y_{p\coprod p'}/C^{p\coprod p'})$$
as in (2.2.1).
Nevertheless (1.8.3) yields an equivalence of
$\Oc^\times$-gerbes
$$C_{p,p'}:
i_{p,p'}^*\bigl(\Detc(Y_{p}/C^{p})\boxtimes\Detc(Y_{p'}/C^{p'})\bigr)\to
\sigma_{p,p'}^*(\Detc(Y_{p\coprod p'}/C^{p\coprod p'}))$$
which is compatible with the
integrable connections along $C^{p\coprod p'}$ in (1.7.3) and takes 
$i_{p,p'}^*(\Ec_p,\Ec_{p'})$ to 
$i_{p,p'}^*(\Ec_p\times\Ec_{p'}),$
viewed as a section of
$$\sigma_{p,p'}^*(\Gr(Y_{p\coprod p'}/C^{p\coprod p'}))=
i_{p,p'}^*(\Gr(Y_p\times Y_{p'}/C^{p}\times C^{p'}))
\to i_{p,p'}^*(Y_p\times Y_{p'})=Y_{p\coprod p'}.$$
Further, applying (1.4.5) to the isomorphism of ind-schemes
$\varkappa_{p,q}$ in (2.2.1) yields an equivalence of 
$\Oc^\times$-gerbes with connections
$$A_{p,q}:
\Delta_{p,q}^*(\Detc_{Y_{pq}/C^{pq}})
\to\Detc_{\Delta_{p,q}^*(Y_{pq})/C^{p}}
\to\varkappa_{p,q}^*(\Detc_{Y_p/C^{p}}).$$
The first arrow is base change while the second one is functoriality, 
see (1.8.4).
In the same way, the isomorphism $\kappa_{p,q}$ yields
an equivalence of $\Oc^\times$-gerbes with connections 
$$B_{p,q}:
j_{p,q}^*(\Detc_{Y_{pq}/C^{pq}})
\to\Detc_{j_{p,q}^*(Y_{pq})/C^{q}}
\to\kappa_{p,q}^*(\Detc_{Y_q/C^{q}}).$$
The compatibility conditions in (2.7.1) are left to the reader.
\qed

\proclaim{(2.7.3) Definition} 
A factorizing object of a
factorization gerbe 
$(\Gc_p)$ 
over the factorization semigroup $(Y_p)$ 
is a collection $(x_p)$ of global objects
equipped with the following data :

(a) for any $p$ in $\Fsetb_0$,
connection data for $x_p$ along $C^{p}$,

(b) for any composable pair 
$p,q$ in $\Fsetb_0$,
isomorphisms 
$\phi_{p,q}: A_{p,q}(\Delta_{p,q}^*(x_{pq}))\to\varkappa_{p,q}^*(x_p)$
and $\varphi_{p,q}: B_{p,q}(j_{p,q}^*(x_{pq}))\to\kappa_{p,q}^*(x_q),$

(c) for any pair $p, p'$ in $\Fsetb_0$, 
isomorphisms $\psi_{p,p'}:C_{p, p'}(i_{p,p'}^*(x_p,x_{p'}))\to
\sigma_{p,p'}^*(x_{p\coprod p'}).$

The data (a)-(c) are required to satisfy the following conditions :

(d) for any composable triple $p,q,r$ in $\Fsetb_0$,
the three diagrams lifting
(2.2.1)(b) are commutative,

(e) for any triple $p, p', p''$ in $\Fsetb_0$, 
the isomorphisms in (c)
are associative and commutative.
\endproclaim

\proclaim{(2.7.4) Proposition} For any two
factorizing objects $(x_p)$, $(y_p)$ of a factorization gerbe $(\Gc_p)$,
the local homomorphisms $x_p\to y_p$ in $\Gc_p$ yield
a system of $\Oc^\times$-torsors $(\Homu (x_p, y_p))$
which is a factorizing line bundle. 
\endproclaim

\noindent{\sl Proof :} 
We supply the data required by Definition 2.6.1. Let $L_p$ be the line
bundle corresponding to the $\Oc^\times$-torsor $(\Homu (x_p, y_p))$. Then:

\vskip .2cm

(a) The connection $\nabla_p$ on $L_p$ comes from the connection data on the objects $x_p$ and $y_p$.

\vskip .1cm

(b) The isomorphism $u_{p,q}: \Delta_{p,q}^*(L_{pq})\to \varkappa_{p,q}^*(L_p)$
comes from the equivalence of gerges $A_{p,q}: \Delta_{p,q}^*(\Gc_{pq})\to \varkappa_{p,q}^*(\Gc_p)$
applied to the Hom-sheaves between $x_{pq}$ and $y_{pq}$. 
Similarly for $v_{p,q}: j_{p,q}^*(L_{pq})\to\kappa_{p,q}^*(L_p)$ 
which come from $B_{p,q}$.

\vskip .1cm

(c) The isomorphisms $w_{p,p'}$ come from the equivalences $C_{p, p'}$. We leave 
further details to the reader. \qed

\vskip 3mm

A factorizing morphism between two factorizing objects $(x_p)$,
$(y_p)$ is a factorizing section of $(\Homu (x_p, y_p))$.
We denote by $(\Gc_p)_{fact}$ the category formed by factorizing
objects and factorizing morphisms of the factorization gerbe $(\Gc_p)$.

\proclaim {(2.7.5) Proposition} 
Assume that $(Y_p)$ is a factorization semigroup of local nature, 
coming from a $(\gen, K)$-ind-scheme $Y$, so that
we have the group of factorizable functions $\Oc^\times(Y)_{fact}$.
Then the category $(\Gc_p)_{fact}$
is an $\Oc^\times(Y)_{fact}$-groupoid, i.e., any Hom-set of this
category is made into an $\Oc^\times(Y)_{fact}$-torsor (possibly
empty) and composition of morphisms is bilinear.
\endproclaim

\vfill\eject

\centerline{\bf 3. The formal loop spaces and chiral differential operators.}

\vskip 1cm

\subhead {(3.1) Reminder on $\Lc^0X$ and $\Lc X$}\endsubhead

For a commutative ring $R$ we denote by 
$R((t))^\surd$ the subring of $R((t))$ consisting of the Laurent series 
$\sum_{i\gg -\infty}^\infty a_i t^i$
such that $a_i$ is nilpotent for $i<0$.

Let $X$ be a scheme of finite type. Recall the scheme $\Lc^0X$ of formal
arcs and the ind-scheme $\Lc X$ of formal loops.
They represent the following functors on $\Schb$ :
$$\Hom_{\Schb}(S, \Lc^0 X) = \Hom_{\Lrsb}\bigl( (S, \Oc_S[[t]]), (X, \Oc_X)\bigr),
$$
$$\Hom_{\Ischb}(S, \Lc X) = 
\Hom_{\Lrsb} \bigl( (S, \Oc_s((t))^\surd), (X, \Oc_X)\bigr).
$$
There are morphisms
$X\buildrel\pi\over\longleftarrow \Lc^0 X\buildrel i\over\longrightarrow
\Lc X,$
where $\pi$ is affine and $i$ 
realizes $\Lc X$ as a formal thickening of $\Lc^0 X$. 

Assume further that $X$ is smooth.
Then $\Lc X$ is smooth locally locally compact. 
If $U\subset X$ is open, then $\Lc U\subset\Lc X$ is open.
Notice that $(\Lc X)_\roman{Zar}=(\Lc^0X)_\roman{Zar}$ since $\Lc X$
is a limit of nilpotent extensions of $\Lc^0X$.
Recall that $\Omega^1_{\Lc X}=\Omega^1_{X,\Lc}$, 
see \cite{KV2, sect.~6.3} for details and notations.
Thus the cotangent sheaf $\Omega^1_{\Lc X}$ 
is locally trivial over $(\Lc X)_\Zar$.
In particular the $\Oc_{\Lc X}^\times$-gerbe $\Detc_{\Lc X}$ is well-defined.

\vskip3mm

\noindent{\bf (3.1.1) Remark.}
Observe that the ind-scheme $\Gr(\Lc X)$ is formally smooth
for any smooth $X$. We'll not use this. Indeed, assume that $X=\AA^N$.
By definition 
$\Gr(\Lc \AA^N)$ is the Sato Grassmaniann
of the trivial vector bundle over $\Lc \AA^N$ with fiber $\CC((t))^N$.
In particular, it is a filtered inductive limit of smooth $\Lc\AA^N$-schemes
(the finite Grassmanians).
Thus it is formally smooth, because a filtered inductive limit
of formally smooth schemes is again formally smooth.

\subhead {(3.2) Reminder on $\Lc_{C^I}^0X$ and $\Lc_{C^I}X$}\endsubhead

Let $C$ be a smooth algebraic curve and $X$ be a smooth algebraic variety.
For any object $I$ of $\Fsetb_0$ we have constructed in
\cite{KV1} a smooth locally locally compact
ind-$C^I$-scheme $\Lc_{C^I}X$ and a $C^I$-subscheme 
$\Lc^0_{C^I}X\subset\Lc_{C^I}X$ such that,
if $I=\{1\}$, the ind-scheme $\Lc_CX$ and the scheme
$\Lc_C^0X$ are obtained by the principal bundle construction
of Gelfand-Kazhdan, i.e.,
$$\Lc_CX=\Lc X\times_K\widehat C,
\quad
\Lc_C^0X=\Lc^0 X\times_K\widehat C.$$

Given a morphism $p:J\to I$ in $\Fsetb_0$, 
we denote by $\Lc_pX$, $\Lc_p^0X$ the restrictions of
$\Lc_{C^J}X$, $\Lc_{C^J}^0X$ to the subscheme $C^{p}$ of $C^J$.
We have the morphisms
$$\Lc^0_pX\buildrel i_p\over\longrightarrow\Lc_pX
\buildrel\rho_p\over\longrightarrow C^{p}.$$
The systems $(\Lc^0_pX)$, $(\Lc_pX)$ 
are factorization semigroups
of local nature 
on $\Lc^0_CX$, $\Lc_CX$ 
by \cite{loc.~cit., (2.3.3), (2.6.2)(b)}.
Further, the
$\Oc_{\Lc_p X}^\times$-gerbe $\Detc_{\Lc_p X/C^p}$ is well-defined
(over the Zariski site) for each $p$. 
Proof is as in the local case in the previous section.

\vskip3mm

\noindent{\bf (3.2.1) Remarks.}
(a)
A factorization monoid is a factorization semigroup $(Y_p)$ with a unit, i.e.,
with a morphism of factorization semigroups $(y_p:C^{p}\to Y_p)$ such that,
for any local section $s$ of $Y_C\to C$, the product
$y_{\{1\}}\times s$ extends to a section of $Y_{\{1,2\}}$ 
(via $\kappa$, $\sigma$)
whose restriction to the diagonal is identified with $s$ 
(via $\varkappa$).
The integrable connection of a factorization monoid can be recovered from 
the other axioms as follows.

Assume that $p=p_{\{1\}}$.
The general case is identical.
Set $I=\{1,2\}$, and let
$\Delta_1,\Delta_2:C^{[I]}\to C$ be the obvious projections.
We must construct an isomorphism of ind-$C^{[I]}$-schemes
$\Delta_1^*(Y_C)\simeq\Delta_2^*(Y_C)$ which restricts to the identity
of $Y_C$ over the diagonal $C\subset C^I$.
By definition of the unit the maps
$\Id\times y_{\{2\}}$, $y_{\{1\}}\times\Id$
yield isomorphisms
$\Delta_1^*(Y_{C}),\Delta_2^*(Y_{C})\to (Y_I)|_{C^{[I]}}$
which restrict to the identity over the diagonal.
This gives a connection as explained in (1.6.2). 
Further, taking $\J = \{1,2,3\}$
and using the unit property gives at once the integrability, see (1.6.2)
as well as Section 3.4.7 of [BD1]. 

(b) 
Although the factorization semigroups $(\Lc^0_pX)$, $(\Lc_pX)$ have no units,
they do have integrable connections along $C^p$.
The integrable connection on $\Lc_pX$ was not given in \cite{KV1}.
For $p=p_{\{1\}}$ this connection follows from the construction
in section 2.3, because
$(\Lc_pX)$ is a factorization semigroup of local nature. 
Let us explain how to get it for an arbitrary $p$.

Once again, to simplify, we assume that $p=p_{\{1\}}$.
The general case is identical.
Set $I=\{1,2\}$.
We must construct an isomorphism of ind-$C^{[I]}$-schemes
$\Delta_1^*(\Lc_{C}X)\to (\Lc_{C^I}X)|_{C^{[I]}}$
which restricts to the identity over the diagonal.
Let $\lambda_{X,C^I}$ be the representating functor of $\Lc_{C^I}X$.
Recall that
$\lambda_{X,C^I}(S)$ consists of pairs $(f_I,\rho)$ such that
$f_I\in Hom(S,C^I)$ and
$\rho\in Hom((\Gamma(f_I),\Kc_{f_I}^\sqr),X)$.
See \cite{KV1, 2.3.2} for details and notations.
An $S$-point of $C^{[I]}$ is a couple $f_I=(f_1,f_2)$ such that
$f_1,f_2:S\to C$ are equal on $S_{red}$.
For such an $f_I$ we must prove that the locally ringed spaces
$(\Gamma(f_I),\Kc_{f_I}^\sqr)$, $(\Gamma(f_1),\Kc_{f_1}^\sqr)$
are isomorphic.
This follows from the following lemma.

\proclaim{(3.2.2) Lemma}
For any $S$-point $f_I$ of $C^{[I]}$ we have
$\widehat\Gamma(f_I)=\widehat\Gamma(f_1)=\widehat\Gamma(f_2)$,
and
$\widehat\Kc_{f_I}=\widehat\Kc_{f_1}=\widehat\Kc_{f_2}$
(as sheaves over 
$\Gamma(f_I)_{red}=\Gamma(f_1)_{red}=\Gamma(f_2)_{red}$).
\endproclaim

\noindent{\sl Proof :}
Choosing an \'etale coordinate on $C$ we reduce to the case $C=\AA^1$,
so $f_1$, $f_2$ are elements of the coordinate ring $B:=\CC[S]$ such that
$s=f_1-f_2$ is a nilpotent element of $B$. Let $n_0$ be such that $s^{n_0}=0$.
Put $R=B[t]=\CC[S\times\AA^1]$ and let $r_i=t-f_i$, so
$r_1-r_2=f_2-f_1=-s$. Then
$$\widehat\Gamma(f_I)=``\ind"\Spec(R/(r_1r_2)^n),\quad
\widehat\Gamma(f_i)=``\ind"\Spec(R/(r_i)^n).$$
On the other hand we have 
$r_1^n\in(r_2^{n-n_0})$
and
$r_2^n\in(r_1^{n-n_0})$.
So the $r_1$-adic and the $r_2$-adic topologies on $R$ are equivalent to
each other and to the $(r_1r_2)$-adic topology.
This implies the first claim.

To prove the second one, let
$$\widehat R=\pro R/(r_1^n)=\pro R/(r_2^n)=\pro R/(r_1r_2)^n.$$
Then $r_1,r_2\in\widehat R$ with $s=r_1-r_2$ nilpotent, 
and our statement means that
$$\widehat R[r_1^{-1}]=
\widehat R[r_2^{-1}]=
\widehat R[(r_1r_2)^{-1}].$$
To see this, let us write
$${1\over r_1}={1\over r_2}\Bigl(1-{s\over r_2}+{s^2\over r_2^2}-\cdots\Bigr)$$
(a terminating geometric series). So $r_1$ is invertible in
$\widehat R[r_2^{-1}]$.
\qed

\vskip3mm

It is easy to check that $\Lc^0_pX$ is a reasonnable subscheme of $\Lc_p X$.
Hence the normal bundle $N_{\Lc_p X/\Lc^0_pX}$ is well-defined. 
It is a vector bundle (of infinite rank) over $\Lc^0_pX$
because $X$ is a smooth variety.
Set $S_{\Lc_pX/\Lc^0_pX}$
equal to the symmetric algebra
$S_{\Oc_{\Lc^0_pX}}(N_{\Lc_pX/\Lc^0_pX})$.
It is proved in \cite{BD1} that
the factorization semigroup $(\Lc^0_pX)$ is commutative.
Further, we have the following.

\proclaim{(3.2.3) Proposition}
The sheaves of $\Oc_{C^p}$-modules
$U\mapsto\rho_{p*}(S_{\Lc_pU/\Lc^0_pU})$
form a sheaf of commutative factorization algebras
over $X$.
\endproclaim

\noindent{\sl Proof :}
We must check that the vertex algebra 
associated to the factorization algebra
$(\rho_{p*}(S_{\Lc_pX/\Lc^0_pX}))$ via (2.4.5) is commutative.
We may assume that $X=\AA^N$.
Then, by \cite{KV1, thm. 5.3.1},
this vertex algebra is isomorphic to $\gr(V_N)$ 
in (3.3.3) below. The claim follows.
\qed

\subhead {(3.3) Reminder on chiral differential operators}\endsubhead

Let $X$ be a smooth scheme of finite type, as before.
So $\Lc^0X$ is a smooth reasonable subscheme of $\Lc X$.
Let $N_{\Lc X/\Lc^0X}$ be the normal bundle of $\Lc^0X$ in $\Lc X$.
It is a vector bundle in $\Ob_{\Lc^0X}$.
The derivation $d/dt$ of $\gen$ makes
$\pi_*(\Oc_{\Lc^0(X)})$ into a sheaf of 
commutative differential algebras on $X$,
and $\pi_*(N_{\Lc X/\Lc^0X})$ into a differential
module over $\pi_*(\Oc_{\Lc^0X})$. So the symmetric algebra
$$S_{\Lc X/\Lc^0X}=\pi_*(S_{\Oc_{\Lc^0X}}(N_{\Lc X/\Lc^0X}))
=S_{\pi_*(\Oc_{\Lc^0X})}(\pi_*(N_{\Lc X/\Lc^0X}))$$
is also a sheaf of commutative differential algebras on $X$.
Therefore, $S_{\Lc X/\Lc^0X}$
can be equipped with a structure of a sheaf of commutative vertex algebras
over $X$ as in (2.4.4).

\proclaim{(3.3.1) Proposition}
The sheaf of factorization algebras 
$U\mapsto(\rho_{p*}(S_{\Lc_pU/\Lc^0_pU}))$ 
gives rise, via (2.4.5), to the
sheaf of vertex algebras
$S_{\Lc X/\Lc^0X}$ above.
\endproclaim

\noindent{\sl Proof :}
Fix an open subset $U\subset X$.
We must prove that the structure of vertex algebra on
$S_{\Lc X/\Lc^0X}(U)$ is the same as the one induced by (3.2.3).
We may assume that $U=\AA^N$.
Then the claim follows from (3.3.3) below.
\qed

\vskip3mm

By a filtration we'll always mean an 
exhaustive positive increasing filtration.

\proclaim{(3.3.2) Definition}
A sheaf of chiral differential operators over $X$ is a sheaf $(\Vc,F)$
of filtered graded vertex algebras, see \cite{KV3, (2.1.3)},
equipped with an isomorphism of sheaves of
graded vertex algebras $gr(\Vc)\to S_{\Lc X/\Lc^0X}.$
A morphism $(\Vc,F)\to (\Vc',F')$ is a morphism of sheaves of filtered 
graded vertex algebras which induces the identity on $S_{\Lc X/\Lc^0X}$.
\endproclaim

In particular we have $F_0(\Vc)=\pi_*(\Oc_{\Lc^0X})$ 
and the isomorphism $gr(\Vc)\to S_{\Lc X/\Lc^0X}$ is identical on $F_0(\Vc)$.
Compare \cite{BD1, (3.9.5)}.

\vskip3mm

\noindent{\bf (3.3.3) Examples.}
(a) Let $X=\AA^N$ and $\Oc^{ch}_{\AA^N}$ be the sheaf of vertex algebras
defined in \cite{MSV}. 
As a Zariski sheaf on $\AA^N$ it is described as follows.
Let $A_N$ be the Heisenberg algebra with generators
$a^i_n, b^i_n,$ $i=1,...N$, $n\in\ZZ$, subject to
the relations $[a^j_m,b^i_n]=\delta_{i,j}\delta_{m,-n}$,
all other brackets being zero.
Consider the cyclic module $V_N=A_N/A_N\{b^i_{<0},a^i_{\le 0}\}$.
Then $\Oc^{ch}_{\AA^N}$ is the
quasicoherent sheaf corresponding to the $\CC[b_0^1,...b_0^n]$-module $V_N$. 
See \cite{MSV} for the vertex structure.
The filtration $F$ on $V_N$ such that
$F_iV_N$ is spanned by the monomials
$$a_{m_1}^{i_1}\cdots a_{m_p}^{i_p}b_{n_1}^{j_1}\cdots b_{n_q}^{j_q},
\quad m_i>0,n_j\ge 0, i\ge p>0, q>0,$$
induces a filtration (also denoted $F$) on
$\Oc^{ch}_{\AA^N}$ and $(\Oc^{ch}_{\AA^N},F)$ is a sheaf of CDO. 
More precisely, $gr(V_N)$ is a commutative vertex algebra and is identified 
with the commutative algebra 
$\CC[a^i_m,b^j_n;m>0,n\ge 0]$
with the derivation
$\partial(a^i_m)=ma^i_{m+1},$
$\partial(b^j_n)=(n+1)b^j_{n+1}.$
See \cite{KV3, (2.1.4)}.

(b) Let $\varphi:U\to\AA^N$ be an \'etale map and
$\Oc^{ch}_{U,\varphi}=\varphi^*(\Oc^{ch}_{\AA^N})$.
Since $\varphi$ induces an isomorphism of formal neighborhoods of points,
methods of \cite{MSV, sect.~3.4} make
$\Oc^{ch}_{U,\varphi}$ into a sheaf of vertex algebras which is a sheaf of CDO
over $U$.

\vskip3mm

Let $\CDOb(X)$ be the groupoid formed by sheaves
of CDO on $X$ and their isomorphisms. The correspondence $U\to\CDOb(U)$ defines 
a stack of groupoids $\CDOc_X$ on the Zariski topology of $X$. 
The main result of \cite{GMS1,2} can be formulated as follows.

\proclaim{(3.3.4) Theorem}
(a) The stack $\CDOc_X$ has lien $\Omega^{2,cl}_X$, the sheaf of closed 2-forms
on $X_\roman{Zar}$. This means that for any two objects
$\Vc,\Wc$ of $\CDOc_X(U)$ the sheaf 
$\Hom_{\CDOc_X(U)}(\Vc,\Wc)$ has a natural structure
of a $\Omega^{2,cl}_U$-torsor (possibly empty) and the composition of
morphisms is $\Omega^{2,cl}_U$-bilinear.

(b) Let $\Pi_0\CDOc_X(U)$ be the set of isomorphism classes of objects of
$\CDOc_X(U)$ and $\Pi_0\CDOc_X$ be the sheaf associated to the presheaf 
$U\mapsto\Pi_0\CDOc_X(U)$.
Then $\Pi_0\CDOc_X$
is identified with $\underline H^3_{DR}$, the sheaf on $X_\roman{Zar}$
corresponding to the presheaf $U\mapsto H^3(\Omega^\bullet(U))$, the third
de Rham cohomology.
\endproclaim

\noindent{\bf (3.3.5) Remark.}
Call a sheaf of CDO over $X$ (Zariski) locally trivial if,
Zariski locally on $X$, it is isomorphic to a sheaf of CDO of the form
$\Oc^{ch}_{U,\varphi}$. 
It follows that the substack $\CDOc_X^{lt}$ of $\CDOc_X$ formed by such
objects is a gerbe with lien $\Omega^{2, cl}_X$.

\subhead{(3.4) Determinantal gerbe and chiral differential operators}
\endsubhead

Let $X$ be a smooth scheme of finite type.
By (2.7.2) the system of gerbes $(\Detc_{\Lc_pX/C^{p}})$ 
has the structure of a factorization 
gerbe over the factorization semigroup $(\Lc_pX)$. 
Further, by (2.7.5) factorizing objects yield
the $\Oc^\times(\Lc X)_{fact}$-groupoid
$(\Detc_{\Lc_pX/C^{p}})_{fact}$.

Let $\pi_*(\Oc^\times_{\Lc X,fact})$ be the sheaf over $X$
associated to the presheaf such that
$U\mapsto\Oc^\times(\Lc U)_{fact}$.
Replacing $X$ by a Zariski open subset $U$ and considering the category 
of factorizing objects, we get a stack 
$(\Detc_{\Lc_pX/C^{p}})_{fact}$
of $\pi_*(\Oc^\times_{\Lc X,fact})$-groupoids over $X$.

\proclaim{(3.4.1) Lemma}
The stack of groupoids 
$(\Detc_{\Lc_pX/C^{p}})_{fact}$
is locally non-empty.
\endproclaim

\noindent{\sl Proof :}
We may assume that $X=\AA^N$ and $C=\AA^1$.
For each $p$ there is a left inverse $\theta_p:\Lc_pX\to\Lc_p^0X$
to the canonical embedding obtained by ``forgetting the negative coefficients
of the formal loops". See \cite{KV1, sect.~2.7} for details.
The inverse image by $\theta_p$ of the relative tangent sheaf of
the smooth compact (ind-) scheme
$\Lc_p^0X$ is a section
of $\Detc_{\Lc_pX/C^{p}}$ over $\Lc_pX.$ 
The collection of all these sections is a factorizing object.
\qed

\vskip3mm

A factorizing local object in
$(\Detc_{\Lc_pX/C^{p}})_{fact}$
is said to be locally trivial if it is locally of the form in the proof
of (3.4.1).
Recall that, for each $p$, the relative tangent sheaf 
of the smooth compact (ind-) scheme
$\Lc_p^0X$ yields a canonical section
of $\Detc_{\Lc_pX/C^{p}}$ over $\Lc_p^0X.$ 
Let $\Detc_{\Lc X,fact}^{lt}$ be the substack of categories of
$(\Detc_{\Lc_pX/C^{p}})_{fact}$
consisting of the locally trivial factorizing sections which restrict 
to the canonical section over $\Lc^0_pX$, and of the factorizing isomorphisms
equal to the identity over $\Lc^0_pX$.
It is again a stack of
$\pi_*(\Oc^\times_{\Lc X,fact})$-groupoids over $X$.
See \cite{KV3, (1.5.4)}.

Set $\o_{p}$ equal to $i_{p\bullet}(\omega_{\Lc_p^0X/C^{p}})$.
For every global object $\Ec_p$ of $\Detc_{\Lc_pX/C^{p}}$ over $\Lc_pX$
the space of sections $\Gamma_{\Ec_p/C^{p}}(\o_{p})$ is a discrete 
$\Oc(\Lc_pX)$-module.
Hence, the direct image by $\rho_p$ is a quasi-coherent sheaf over $C^{p}$.
Replacing $X$ by open subsets, it localizes to a sheaf
$\underline\Gamma_{\Ec_p/C^{p}}(\o_{p})$
over $X\times C^p$.

Finally, recall that
$\CDOc_X^{lt}$ is a gerbe on $X$ with lien equal to the sheaf
$\Omega^{2, cl}_X$.
In \cite{KV3, (2.4.3)}
we construct an identification of 
sheaves of Abelian groups
$S:\Omega^{2,cl}_X\to\pi_*(\Oc^\times_{\Lc X,fact})$
via the multiplicative symplectic action map.

We can now prove the main result of this paper.

\proclaim{(3.4.2) Theorem}
(a) 
For each global object 
$\Ec=\Ec_C,(\Ec_p)$ of 
$\Detc_{\Lc X,fact}^{lt}$
the system of sheaves 
$(\underline\Gamma_{\Ec_p/C^{p}}(\o_{p}))$ 
is a sheaf of factorization algebras over $X$.

(b) 
Fix a closed point $x\in C$.
The fiber of
$\underline\Gamma_{\Ec_{\{1\}}/C}(\o_{\{1\}})$ 
at $x$ is a sheaf of vertex algebras over $X$.
It is equipped with a natural filtration making it into a sheaf
of CDO over $X$. 

(c) The functor $\Ec\mapsto\Gamma_\Ec$ is an equivalence
of gerbes $\Detc_{\Lc X,fact}^{lt}\to\CDOc_X^{lt}$
which is compatible with the identification of their liens by $S^{-1}$.
\endproclaim

\noindent{\sl Proof :}
Replacing $X$ by open subsets, we must construct isomorphisms
$$\aligned
&\mu_{p,q}:\Delta_{p,q}^*(\rho_{pq*}\Gamma_{\Ec_{pq}/C^{pq}}(\o_{pq}))
\to\rho_{p*}\Gamma_{\Ec_p/C^{p}}(\o_{p}),
\\
&\l_{p,q}:j_{p,q}^*(\rho_{pq*}\Gamma_{\Ec_{pq}/C^{pq}}(\o_{pq}))
\to\rho_{q*}\Gamma_{\Ec_q/C^{q}}(\o_{q}),
\\
&\nu_{p,p'}:
i_{p,p'}^*\bigl(
\rho_{p*}\Gamma_{\Ec_p/C^{p}}(\o_{p})\boxtimes
\rho_{p'*}\Gamma_{\Ec_{p'}/C^{p'}}(\o_{p'})\bigr)\to
\\
&\qquad\qquad\qquad\to
\rho_{p\coprod p'*}\Gamma_{\Ec_{p\coprod p'}/C^{p\coprod p'}}
(\o_{p\coprod p'})
\endaligned$$
and a global section (=the unit) of
$\Gamma_{\Ec_{\{1\}}/C}(\o_{\{1\}})$ 
satisfying the axioms in (2.4.1).

The factorization semigroup on $\Lc X$
yields isomorphisms of ind-schemes
$$\varkappa_{p,q}:\Delta_{p,q}^*(\Lc_{pq}X)\to\Lc_pX,
\quad
\kappa_{p,q}:j_{p,q}^*(\Lc_{pq}X)\to\Lc_qX,$$
$$\sigma_{p,p'}:i_{p,p'}^*(\Lc_pX\times\Lc_{p'}X)\to\Lc_{p\coprod p'}X.$$
We'll use the same notations for the scheme $\Lc^0X$.
By \cite{KV1, proof of (2.4.1)} the scheme 
$\Lc_p^0X$ is the limit of a filtering projective system
$(\Lc^0_{p,n}X)$ consisting of smooth schemes of finite type. 
Further, the relative $\Dc$-module 
$\omega_{\Lc_p^0X/C^{p}}$ is represented by the 
relative $\Dc$-module 
$\omega_{\Lc_{p,n}^0X/C^{p}}$ on $\Lc_{p,n}^0X$ for any $n$.
So the relative $\Dc$-modules
$\sigma_{p,p'}^*(\o_{p\coprod p'})$ and 
$i^*_{p,p'}(\o_{p}\boxtimes\o_{p'})$
are isomorphic.

Recall that $\Ec_p$ is a section of the relative Sato Grassmannian
$\Gr(\Lc_pX/C^{p})$ over $\Lc_pX$
and that the equivalence $C_{p,p'}$
takes $i^*_{p,p'}(\Ec_p,\Ec_{p'})$ to the section
$i^*_{p,p'}(\Ec_p\times\Ec_{p'})$ of 
$\sigma_{p,p'}^*(\Gr(\Lc_{p\coprod p'}X/C^{p\coprod p'}))$.
As $(\Ec_p)$ is a factorizing object,
the later is identified with $\sigma_{p,p'}^*(\Ec_{p\coprod p'})$ 
by $\psi_{p,p'}$.
Therefore (1.8.3) yields an isomorphism
$$
\nu_{p,p'}:
i_{p,p'}^*\bigl(
\Gamma_{\Ec_p/C^{p}}(\o_{p})\boxtimes
\Gamma_{\Ec_{p'}/C^{p'}}(\o_{p'})\bigr)
\to$$
$$\to
i_{p,p'}^*\bigl(
\Gamma_{\Ec_p\times\Ec_{p'}/C^{p}\times C^{p'}}
(\o_{p}\boxtimes\o_{p'})\bigr)
\to
\Gamma_{\Ec_{p\coprod p'}/C^{p\coprod p'}}(\o_{p\coprod p'}).
$$
The construction of the isomorphisms $\mu_{p,q}$, $\l_{p,q}$ 
and the proof of the compatibility conditions
are left to the reader.
The unit is the delta function $\delta_{\Lc^0_CX}$ in 
$\Gamma_{\Ec_{\{1\}}/C}(\o_{\{1\}})$, see below for more details.

Now we concentrate on claim (b).
It is enough to construct a filtration of 
$\Gamma_{\Ec_p/C^{p}}(\o_{p})$
which is compatible with the factorization structure 
for $X$ affine with an \'etale map $\phi:X\to\AA^N$ and $C=\AA^1$.
By \cite{loc.~ cit., (2.9.3)} the ind-scheme $\Lc_p X$ is represented 
in the following way
$$\Lc_p X={``\ind"}_{\!\!\!\eps}\pro_n\Lc_{p,n}^\eps\phi
=\pro_n\Lc_{p,n}^\infty\phi,$$
where $\Lc_{p,n}^\eps \phi$ is a scheme of finite type
and 
$\Lc_{p,n}^\infty\phi={``\ind"}_{\!\!\!\eps}\Lc_{p,n}^\eps \phi.$
To simplify we may further assume that
$X=\AA^d$ and $\phi=\Id$.
Then, 
by \cite{loc.~ cit., (5.5.5)},
there are closed embeddings
$\Lc^0_{p,n}X\subset\Lc^\eps_{p,n}\phi\subset X^\eps_{p,n}$,
where $X^\eps_{p,n}$ is an affine space of finite dimension.
Let $i^{0\eps}$ be the first inclusion,
and 
$i^{0X}$, $i^{\eps X}$ be the inclusions 
$\Lc^0_{p,n}X,$ $\Lc^\eps_{p,n}\phi\subset X^\eps_{p,n}$.

Consider the right $\Dc_{\Lc^\eps_{p,n}\phi/C^{p}}$-module
$\o_{p,n,\eps}=i_{\bullet}^{0\eps}(\o_{\Lc^0_{p,n}X/C^{p}})$.
As $\Lc^\eps_{p,n}\phi$ is a singular scheme,
$\o_{p,n,\eps}$ must be viewed as a
$\Dc_{X^\eps_{p,n}/C^{p}}$-module supported on $\Lc^\eps_{p,n}\phi,$
i.e., it is equal to
$(i^{\eps X})^!i_{\bullet}^{0X}(\o_{\Lc^0_{p,n}X/C^{p}})$.

The limit $\Lc^\infty_{p,n}\phi$ of the filtering inductive system of schemes
$(\Lc^\eps_{p,n}\phi)$ is a formally smooth ind-scheme of ind-finite type.
Let $\pi_{nn'}$ be the natural projection
$\Lc^\infty_{p,n'}\phi\to\Lc^\infty_{p,n}\phi$.
Then $\Ec_p$ consists of a family of subbundles of finite rank
$\Ec_{p,n}$ of $\Theta_{\Lc^\infty_{p,n}\phi/C^{p}}$ 
such that $\Ec_{p,n'}$ is equal to
$(d\pi_{nn'})^{-1}(\pi_{nn'}^*(\Ec_{p,n}))$
for each $n'\ge n$.
Further the space of sections 
$\Gamma_{\Ec_p/C^{p}}(\o_{p})$
is the limit of the filtering inductive system of vector spaces
$(\Gamma_{p,n}^\eps)$
such that
$$\Gamma_{p,n}^\eps=
\Gamma\Bigl(\Lc^{\eps}_{p,n}\phi,
\o_{p,n,\eps}
\otimes_{\Oc_{\Lc^{\eps}_{p,n}\phi}}
\Wedge^{\max}(\Ec_{p,n})|_{\Lc^{\eps}_{p,n}\phi}\Bigr).$$ 
Let $D_i$ be the subspace of 
$\Dc_{X^\eps_{p,n}/C^{p}}(X^\eps_{p,n})$
consisting of the differential operators of degree $\le i$,
and set $\o_i$ equal to the subsheaf
$$(i_*^{0X}(\o_{\Lc^0_{p,n}X/C^{p}}))\cdot D_i\subset
i_\bullet^{0X}(\o_{\Lc^0_{p,n}X/C^{p}}).$$
Let $F_i(\Gamma_{p,n}^\eps)$ be the subspace of $\Gamma_{p,n}^\eps$ equal to 
$$\Gamma\Bigl(\Lc^{\eps}_{p,n}\phi,
(i^{\eps X})^!(\o_i)
\otimes_{\Oc_{\Lc^{\eps}_{p,n}\phi}}
\Wedge^{\max}(\Ec_{p,n})|_{\Lc^{\eps}_{p,n}\phi}\Bigr)
.$$ 
The same argument as in (1.8) shows that the vector spaces
$F_i(\Gamma_{p,n}^\eps)$ form again
a filtering inductive system 
for each $i,p$.
Thus we may set
$$F_i(\Gamma_{\Ec_p/C^{p}}(\o_{p}))=
\ind_{\eps,n}F_i(\Gamma_{p,n}^\eps).$$
This is the required filtration
on $\Gamma_{\Ec_p/C^{p}}(\o_{p})$.
The compatibility with the factorization structure is proved as in (a) above.
Setting $p=p_{\{1\}}$ this filtration yields a filtration on
$\Gamma_\Ec$ in the obvious way.

We can now define the unit
of the factorization algebra $\Gamma_{\Ec_p/C^{p}}(\o_{p})$
mentioned in (a) above.
Setting $p$ equal to the map $p_{\{1\}}$,
we get the scheme $\Lc^0_{C,0}X$ which is equal to $X\times C$ and
the sheaf $\o_{\Lc,C,0,0}$ which is equal to $\o_X\boxtimes\Oc_C$.
By definition of $\Ec$ we have
$\Ec_{p,n}|_{\Lc^0_{p,n}X}=\Theta_{\Lc^0_{p,n}X/C^p}.$ 
In particular we have  $\Ec_{C,0}|_{X\times C}=\Theta_X.$ 
Thus the space $\Gamma_{p,0}^0$ is equal to 
$\Gamma(X\times C,\Oc_{X\times C})$ and the delta function
$\delta_{\Lc^0_CX}$ in
$\Gamma_{\Ec_{\{1\}}/C}(\o_{\{1\}})$
is identified with the unit in 
$\Gamma(X\times C,\Oc_{X\times C})$.

Next, by (3.3.1) we must prove that the factorization algebras
$(\rho_{p*}\gr(\Gamma_{\Ec_p/C^p}(\o_{p})))$
and $(\rho_{p*}(S_{\Lc_p X/\Lc^0_pX}))$ are isomorphic.
For each integer $n\ge 0$ set 
$F_i(\Gamma_{p,n}^\infty)=\ind F_i(\Gamma_{p,n}^\eps).$
The graded vector space 
$$\Oplus_iF_{i+1}(\Gamma_{p,n}^\infty)/F_i(\Gamma_{p,n}^\infty)$$
is the symmetric algebra
of the $\CC[\Lc^0_{p,n}X]$-module
$$N_{\Lc^\infty_{p,n}\phi/\Lc^0_{p,n}X}=
(\Theta_{\Lc^\infty_{p,n}\phi}|_{\Lc^0_{p,n}X})/\Theta_{\Lc^0_{p,n}X}.$$
On the other hand 
$N_{\Lc_p X/\Lc^0_pX}$ is the pull-back of the bundle
$N_{\Lc^\infty_{p,n}\phi/\Lc^0_{p,n}X}$
by the natural projection $\Lc^0_p X\to\Lc^0_{p,n}X$.
So
$\gr(\Gamma_{\Ec_p/C^{p}}(\o_{p}))$
is equal to $S_{\Lc_p X/\Lc^0_pX}.$
The compatibility with the factorization structures is routine.

Finally, let us concentrate on claim (c).
Functoriality of $\Ec\mapsto\Gamma_\Ec$ follows from (1.8.2).
We have also proved there that the local automorphism
$\Ec\to\Ec$
associated to a local section 
$f$ of $\pi_*(\Oc^\times_{\Lc X,fact})$
is taken to the automorphism 
$\Gamma_\Ec\to\Gamma_\Ec$
given by multiplication by $f$.
Since the latter is equal to 
the action 
in \cite{GMS1} 
of the closed 2-form $S^{-1}(f)$ 
on the sheaf of chiral differential operators 
$\Gamma_\Ec$ over $X$,
see \cite{KV3, (2.3.4)},
the functor 
$\Gamma:\Detc_{\Lc X,fact}^{lt}\to\CDOc_X$
is compatible with the identification 
of the liens by $S^{-1}$.
We are done, because $\Gamma$ maps into the gerbe $\CDOc_X^{lt}.$
\qed

\vskip3mm

In \cite{KV2, sect.~4.3, 6.2} 
we construct
a pull-back homomorphism
$$ev^*:H^2\bigl(X,K_2(\Oc_X)\bigr)\to 
H^2\bigl(\Lc X,K_2(\Oc_{\Lc X}((t)))\bigr)$$
by the morphism of ringed spaces
$ev:(\Lc^0X,\Oc_{\Lc X}((t)))\to(X,\Oc_X)$
as well as the Contou-Carr\`ere symbol
$$\partial : H^2\bigl(\Lc X,K_2 (\Oc_{\Lc X}((t)))\bigr)\to
H^2\bigl(\Lc X,\Oc_{\Lc X}^\times\bigr)=
H^2\bigl(X,\pi_*(\Oc_{\Lc X}^\times)\bigr).$$
Let $ch_2(\Theta_X)$ be the Chern character in 
$H^2(X,K_2(\Oc_X))$, and $[\CDOc_X^{lt}]$
be the class in $H^2(X,\Omega^{2, cl}_X)$.

\proclaim{(3.4.3) Corollary}
The classes
$S_*[\CDOc_X^{lt}]$ and $\partial(ev^*(ch_2(\Theta_X)))$ 
in $H^2(X,\pi_*(\Oc^\times_{\Lc X}))\otimes\ZZ[{1\over 2}]$
are equal.
\endproclaim

\noindent{\sl Proof :}
The classes $S_*[\CDOc_X^{lt}]$ 
and $[\Detc_{\Lc X,fact}^{lt}]$ 
in $H^2(X,\pi_*(\Oc^\times_{\Lc X,fact}))$ are equal.
The embedding of sheaves
$\pi_*(\Oc^\times_{\Lc X,fact})\subset\pi_*(\Oc^\times_{\Lc X})$
yields a map
$$H^2(X,\pi_*(\Oc^\times_{\Lc X,fact}))\to H^2(X,\pi_*(\Oc^\times_{\Lc X})).$$
Since the $\pi_*(\Oc^\times_{\Lc X})$-gerbes
$\pi_*(\Detc_{\Lc X})$ and 
$$\Detc_{\Lc X,fact}^{lt}\otimes_
{\pi_*(\Oc^\times_{\Lc X,fact})}\pi_*(\Oc^\times_{\Lc X})$$
are isomorphic,
by \cite{KV2, (6.3.3), (6.4.1)} 
the image of
$[\Detc_{\Lc X,fact}^{lt}]$ 
in $H^2(X,\pi_*(\Oc^\times_{\Lc X}))\otimes\ZZ[{1\over 2}]$ 
is equal to
$\partial(ev^*(ch_2(\Theta_X)))$.
Therefore we get the equality of classes
$$S_*[\CDOc_X^{lt}]=\partial(ev^*(ch_2(\Theta_X))).$$
\qed


\Refs
\widestnumber\key{ABCD}

\ref\key{BD1}\by Beilinson, A., Drinfeld, V.
\book Chiral Algebras 
\bookinfo \vol 
\publ American Mathematical Society\yr 2004\endref

\ref\key{BD2}\by Beilinson, A., Drinfeld, V.
\book Quantization of Hitchin's Hamiltonians and Hecke eigensheaves 
\bookinfo preprint\vol 
\publ \yr \endref

\ref\key{BM}\by Breen, L., Messing, W. 
\paper Differential geometry of gerbes   
\vol 198
\jour Adv. in Math. 
\yr 2005
\pages 732--846
\endref



\ref\key{D}\by Drinfeld, V.
\paper Infinite-dimensional vector bundles in algebraic geometry 
(an introduction) 
 \book The unity of mathematics 
\bookinfo  Progr. Math., 244
\publ Birkhäuser Boston
\yr 2006\endref


 

\ref\key[E] \by Emery, M. \book Stochastic Calculus on Manifolds
\publ. Springer-Verlag \yr 1989\endref

\ref\key{GMS1}\by Gorbounov, V., Malikov, F., Schechtman, V.
\paper Gerbes of chiral differential operators 
\jour Mathematical Research Letters
\vol 7
\yr 2000
\pages 55-66
\endref

\ref\key{GMS2}\by Gorbounov, V., Malikov, F., Schechtman, V.
\paper Gerbes of chiral differential operators. II. Vertex algebroids
\jour Invent. Math.
\vol  155  
\yr 2004
\pages 605--680 
\endref

\ref\key{K}\by Kac, V.\book Vertex algebras for beginers
\bookinfo University Lecture Series\vol 10
\publ American Mathematical Society\yr 1996\endref

\ref\key{KV1}\by Kapranov, M., Vasserot, E.
\paper Vertex algebras and the formal loop space
\jour Publ. Math., Inst. Hautes Etud. Sci.
\vol 100
\yr 2004
\pages  209--269.
\endref

\ref\key{KV2}\by Kapranov, M., Vasserot, E.
\paper Formal Loops II : the local Riemann-Roch theorem
for determinantal gerbes
\jour Ann. Sci. ENS, to appear
\vol
\yr
\endref

\ref\key{KV3}\by Kapranov, M., Vasserot, E.
\paper Formal Loops III : factorizing functions and the Radon transform 
\jour math.AG/0510476
\vol
\yr
\endref

\ref\key{McK} \by McKean, H. P. \book Stochastic Integrals
\publ Chelsea Publ. \yr 2005\endref

\ref\key{M}\by Mac Lane, S. 
\book Categories For the Working Mathematician. Second edition
\bookinfo GTM\vol 5
\publ Springer\yr 1998\endref

\ref\key{MSV}\by  Malikov, F., Schechtman, V., Vaintrob, A.
\paper Chiral
de Rham complex
\jour Comm. Math. Phys. 
\vol  204  
\yr 1999 
\pages 439-473
\endref

\ref\key{PS}\by Pressley, A., Segal, G.B.\book Loop Groups
\publ Cambridge Univ. Press\yr 1986\endref
\bye